\documentclass[12pt]{article}

\usepackage[raggedright]{subfigure}

\usepackage{geometry}
\usepackage{amsmath}
\usepackage{amssymb}
\usepackage{amsthm}

\usepackage{eepic}

\newcommand{\dist}{d}

\newtheorem{thm}{Theorem}[section] 
\newtheorem{thmi}{Theorem}[section]
\newtheorem{cor}[thm]{Corollary} 
\newtheorem*{cor2}{Corollary to Lemma~\ref{wmalemma1}}
\newtheorem{lm}{Lemma}[section]

\theoremstyle{definition}
\newtheorem{define}{Definition}[section] 
\newtheorem{cl}{Claim}[lm]

\theoremstyle{remark}
\newtheorem*{rem}{Remark}

\begin{document}

\title{The DNA Inequality in Non-Convex Regions}
\author{Eric Larson}
\date{}
\maketitle

\begin{abstract}
A simple plane closed curve $\Gamma$ satisfies the \emph{DNA Inequality} if the
average curvature of any closed curve contained inside $\Gamma$ exceeds
the average curvature of $\Gamma$. In 1997 Lagarias and Richardson
proved that all convex curves satisfy the DNA Inequality and asked
whether this is true for some non-convex curve. They conjectured that
the DNA Inequality holds for certain L-shaped curves. In this paper,
we disprove this conjecture for all L-Shapes and construct a large class
of non-convex curves for which the DNA Inequality holds. We also give
a polynomial-time procedure for determining whether any specific curve
in a much larger class satisfies the DNA Inequality.
\end{abstract}

\section{\label{secintro} Introduction}

A simple plane closed curve $\Gamma$ is said to satisfy the \emph{DNA Inequality}
if the average curvature (which is the integral of the absolute value of
curvature divided by the perimeter) of any closed curve contained within
the region bounded by $\Gamma$ exceeds the average curvature of $\Gamma$.
(It is called the ``DNA Inequality" because the picture is akin to a
little piece of DNA inside of a cell.) In the following, we will refer
to the outside closed curve $\Gamma$ as the ``cell,'' and the inside
closed curve as the ``DNA'' (denoted $\gamma$). All cells considered in
this paper will be (non-self-intersecting) closed polygons, but the DNA
closed curves are allowed to have self-intersections. The DNA Inequality
has been proven to hold for all convex cells; see \cite{L-R, N-P,T}.
On the second page of the paper by Lagarias and Richardson \cite{L-R}
that proved it for convex cells, they raised the question whether
the DNA Inequality might hold for some non-convex cells. In particular,
they suggested that some L-shaped regions might satisfy the DNA
inequality. Here an L-shaped region is a rectangle with a smaller
rectangle removed from one corner of it, cf.\ Section~\ref{seclshape}.

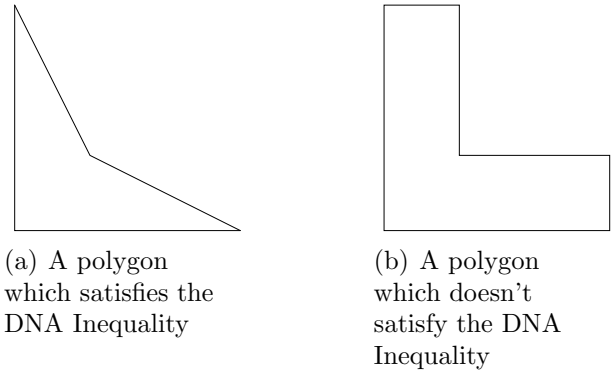
\begin{figure}[ht] \label{figex}
\begin{center}
\subfigure[A polygon which satisfies the DNA Inequality]{\label{sdent}
\unitlength 1.0 mm
\begin{picture}(30,30)(0,0)
\path(30,0)(0,0)(0,30)
\path(0,30)(10,10)(30,0)
\end{picture}}
\hspace{1.5cm}
\subfigure[A polygon which doesn't satisfy the DNA Inequality]{\label{slsh}
\unitlength 1.0 mm
\begin{picture}(30,30)(0,0)
\path(0,0)(0,30)(10,30)(10,10)(30,10)(30,0)(0,0)
\end{picture}}
\end{center}
\caption{Some Polygons}
\end{figure}

The question of when the DNA inequality might hold for non-convex cells is the 
focus of this paper. We obtain three main results, stated below. In particular, 
for the polygons pictured in Figure~\ref{figex}, our results imply that the 
L-shaped polygon \ref{slsh} does not satisfy the DNA inequality, however the non-convex 
quadrilateral \ref{sdent} does satisfy the DNA inequality.

\smallskip

\noindent
Our first result is as follows.
\begin{thmi} \label{ilshape} (Theorem~\ref{lshape}.)
The DNA inequality is false for all L-Shapes.
\end{thmi}

This result disproves the suggestion of Lagarias and Richarson [Conjecture~p.2,~1]. 
This is shown in Section~\ref{seclshape}.

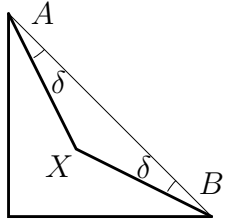
\begin{figure}[ht]
\begin{center}
\unitlength 0.9 mm
\begin{picture}(35,37.5)(0,0)
\Thicklines
\path(30,0)(0,0)(0,30)
\path(0,30)(10,10)(30,0)
\thinlines
\path(0,30)(30,0)
\put(0,30){\arc{15}{0.7854}{1.0472}}
\put(30,0){\arc{15}{3.6652}{3.927}}
\put(7.4,20){\makebox(0,0)[cc]{$\delta$}}
\put(7.5,7.5){\makebox(0,0)[cc]{$X$}}
\put(20,7.4){\makebox(0,0)[cc]{$\delta$}}
\put(5,30){\makebox(0,0)[cc]{$A$}}
\put(30,5){\makebox(0,0)[cc]{$B$}}
\end{picture}
\end{center}
\caption{Example of $P_\delta$} \label{figPdelta}
\end{figure}

Our second result is the main result of this paper. It shows that the DNA
Inequality does hold for a class of non-convex polygonal cells. These are
cells obtained from particular
convex polygons by putting an ``isosceles dent'' in a particular one of
its sides. Namely, take some convex polygon $P$, and fix a side $AB$ of
that polygon. Construct a point $X$ such that 
$\angle{XAB} = \angle{XBA} = \delta$. 
In Figure~\ref{figPdelta} is pictured this construction when $P$ is an
isosceles right triangle, and $AB$ is the hypotenuse.

\begin{define} \label{Pdelta} 
For any convex polygon $P$ with a fixed side $AB$, we denote the curve
which is created from $P$, replacing $AB$ with the two segments $AX$
and $XB$, which is pictured in Figure~\ref{figPdelta} as the bold curve,
by $P_\delta$. (Thus, $P_\delta$ is an ``isosceles denting of $P$ along
$AB$.'') Given a convex polygon $P$ and an edge $E$, we call $(P,E)$ a
\emph{deformable DNA-polygon} (called a \emph{DDNA-polygon} for short)
if there exists a $\delta_0 > 0$ such that $\delta \leq \delta_0$
implies that $P_\delta$ satisfies the DNA Inequality.
\end{define}

\noindent
In this paper, we both classify all DDNA-polygons and create a polynomial
time algorithm for determining if a given curve in a larger class
satisfies the DNA Inequality.

\medskip

\noindent
The following result classifies all DDNA-polygons.

\begin{thmi} \label{ismalldent} (Theorem~\ref{smalldent}.)
If $P$ is a convex polygon with perimeter $p$ and we are denting an
edge with length $l$, and $\alpha$ is the larger of the two angles that
the edge makes with the two adjacent edges, then $P$ is a DDNA-polygon
(with respect to this edge) if and only if:
$$2p \leq \pi l \frac{1 + \cos{\alpha}}{\sin{\alpha}}.$$
\end{thmi}

\begin{rem}
In this paper, we assume for simplicity that the dent $XAB$ is isosceles. The
methods of this paper can still be applied if the triangle is not
isosceles, or even if there are multiple dents all depending on one
parameter $\delta$, so long as adjacent sides are not dented and there
exists $\epsilon > 0$, which does not depend on $\delta$, such that any
two vertices of $P_\delta$ are at least $\epsilon$ apart.
\end{rem}

Our third result is algorithmic, and applies to a class of non-convex polygons 
which we term separable polygons. To define these, we say that an
\emph{interior vertex} of a polygon $\Gamma$ is a vertex contained in
the interior of the convex hull of $\Gamma$. A \emph{separable polygon}
$\Gamma$ is a polygon having the property that for any point $p$ in
the interior of the cell determined by $\Gamma$, but not a vertex of
$\Gamma$, there is at most one interior vertex $v$ of $\Gamma$ such
that the straight line determined by $p$ and $v$ intersects $\Gamma$ in
more than two points. We find a polynomial time algorithm, which when
given a (non-convex) separable polygon determines whether it satisfies the
DNA Inequality.

\begin{thmi} \label{ipoly} (Theorem~\ref{poly}.)
There exists an algorithm which when given as input a separable 
polygon $\Gamma$ specified by its $n$ vertices,
determines whether or not $\Gamma$ satisfies the DNA 
inequality, in number of elementary operations polynomial in $n$.
\end{thmi}

We analyze this algorithm using a simplified model of computation described in 
Section~\ref{secpoly}, in which an ``elementary operation'' is defined.

In Section~\ref{seclshape}, we prove Theorem~\ref{ilshape}. In
Section~\ref{secnot}, we set up the notation that we will use for the
proof of Theorem~\ref{ismalldent}, and give an outline of the proof. In
Section~\ref{seclems}, we prove some useful Lemmas that apply to any
cell. In Section~\ref{secsep}, we turn our attention to a special
class of polygons (which we term ``separable polygons''); we prove
that the DNA Inequality holds in any separable polygon if and only if
it holds for some specific types of DNA. In Section~\ref{secpoly},
we see that this produces a polynomial-time algorithm to determine
whether any separable polygon satisfies the DNA Inequality
(Theorem~\ref{ipoly}). In Section~\ref{secseq}, we determine, using the
results from Section~\ref{secsep} as well as the results of \cite{L-R} and
\cite{N-P}, what happens when we have a sequence of non-convex polygons
which approach a convex one. Finally, in Section~\ref{secsmalldent} we
state and prove Theorem~\ref{ismalldent}, and a corollary which gives some
sufficient and necessary conditions for a polygon to be a DDNA-polygon.

\paragraph{Acknowledgements.}
This research started at the Penn State REU, supported by NSF Grant
No.\ 0505430. I would like to thank Misha Guysinsky and Serge
Tabachnikov for bringing the this problem to my attention, Ken Ross
for helpful discussions and help editing this paper, and the
anonymous referees for their useful critique and numerous suggestions
for improving the presentation.

\section{\label{seclshape} Disproof of the DNA Inequality for L-Shapes}

Here, we present the proof of Theorem~\ref{ilshape}.

\begin{thm}\label{lshape} 
The DNA inequality is false for all L-Shapes.
\end{thm}

\begin{proof}

\looseness=-1
We proceed to construct a counterexample to the DNA Inequality for any
L-Shape.  Choose some sufficiently small $\theta$.  (The size of $\theta$
is bounded above by the dimensions of the L-Shape, but it will be clear
that some nonzero $\theta$ can always be chosen.)  Construct points 
$P \in (A, B)$ and $Q \in (C, D)$ such that $\angle AYP = \angle DYQ = \theta$.

\begin{figure}[ht]
\begin{center}
\unitlength 0.9 mm
\begin{picture}(60,60)(0,0)
\path(0,0)(0,60)(60,60)(60,0)(0,0)
\path(20,60)(20,20)(60,20)
\path(0,0)(10,60)
\path(0,0)(60,10)
\put(62,62){\makebox(0,0)[cc]{$Z$}}
\put(0,62){\makebox(0,0)[cc]{$A$}}
\put(10,62){\makebox(0,0)[cc]{$P$}}
\put(20,62){\makebox(0,0)[cc]{$B$}}
\put(62,10){\makebox(0,0)[cc]{$Q$}}
\put(62,20){\makebox(0,0)[cc]{$C$}}
\put(62,0){\makebox(0,0)[cc]{$D$}}
\put(22,22){\makebox(0,0)[cc]{$X$}}
\put(3,3){\makebox(0,0)[cc]{$Y$}}
\put(2,28){\makebox(0,0)[cc]{$\theta$}}
\put(28,2){\makebox(0,0)[cc]{$\theta$}}
\end{picture}
\end{center}
\caption{Counterexample for L-Shapes} \label{figcounter}
\end{figure}
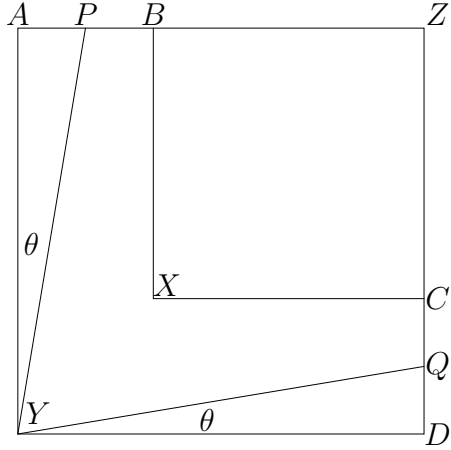

We consider the closed curve $A, P, Y, Q, D, Y, A$. (see
Figure~\ref{figcounter}) Its curvature is clearly $3\pi + 4\theta$, and
its perimeter is clearly $(AY + YD)(1 + \sec(\theta) + \tan(\theta))$.
The curvature of the whole figure is $3\pi$, and the perimeter is 
$2(AY + YD)$. Therefore, to disprove the DNA Inequality, we will show:
\begin{align*}
\frac{3\pi + 4\theta}{(AY + YD)(1 + \sec(\theta) + \tan(\theta))} &< \frac{3\pi}{2(AY+YD)} \\
\Longleftrightarrow \frac{3\pi + 4\theta}{1 + \sec(\theta) + \tan(\theta)} &< \frac{3\pi}{2}\\
\Longleftarrow \frac{3\pi + 4\theta}{2 + \tan(\theta)} &< \frac{3\pi}{2}\\
\Longleftrightarrow \frac{8}{3\pi} &<  \frac{\tan(\theta)}{\theta}
\end{align*}
To verify this, it suffices to note that:
$$\frac{8}{3\pi} < 1 \leq \frac{\tan(\theta)}{\theta}$$
Thus, the DNA Inequality is false for all L-Shapes.

\end{proof}

\begin{rem}
Even if one were to require that the DNA was not self-intersecting, one
could still construct a counterexample by moving the vertex of the curve
that we constructed above coinciding with $Y$, which occurs between $P$
and $Q$, a tiny bit towards $X$.
\end{rem}

\section{\label{secnot} Outline of the Proof of Theorem~\ref{ismalldent}}

The proof of Theorem~\ref{ismalldent} is quite involved, 
and in the course of proving it we establish Theorem~\ref{ipoly}. We first set up some 
basic notation and then outline the proof of Theorem~\ref{ismalldent}.

To prove the DNA Inequality for any curve, it suffices to prove it for
closed polygonal lines. In this case, the integral of the absolute
curvature reduces to a sum of the exterior angles at the vertices (where
the exterior angles are measured so that they are in the interval
$[0,\pi]$). For an explanation of this reduction see \cite{L-R},
Section~2. (In particular, equation (2.3).)

We set our notation for polygonal curves.
We write $\gamma$ for the closed polygonal ``DNA.'' We
denote the vertex sequence of $\gamma$ by 
$\gamma_0, \gamma_1, \ldots, \gamma_n = \gamma_0$, 
and the vertex sequence of $\Gamma$ by 
$\Gamma_0, \Gamma_1, \ldots, \Gamma_m = \Gamma_0$.
When we refer to the number of vertices of some polygon, we shall
mean the number of vertices with multiplicity, unless otherwise stated.
We consider indices modulo $n$ (modulo $m$ for $\Gamma$), and assume that
we never have $\gamma_i, \gamma_{i+1}, \gamma_{i+2}$ collinear. (Under
this assumption, the exterior angles are in the interval $(0, \pi]$.)

\begin{define} \label{f} We define:
$$f_\Gamma(\gamma) = \alpha \cdot (\mbox{curvature of $\gamma$}) - (\mbox{perimeter of $\gamma$})$$
where $1/\alpha$ is the average curvature of $\Gamma$. 
\end{define}

Of course, $\Gamma$ satisfies the DNA Inequality means that
$f_\Gamma(\gamma) \geq 0$ for any closed curve $\gamma$ contained in
$\Gamma$.

\begin{define} \label{CX} 
We term a closed polygonal DNA $\gamma$ contained within the cell
$\Gamma$ with $f_\Gamma(\gamma) < 0$ a \emph{$CX_\Gamma$-polygon.}
(As it is a ``counterexample'' to the DNA Inequality in $\Gamma$.) 
\end{define}

\begin{define} \label{ell} 
We write $\dist(X,Y)$ for the distance between $X$ and $Y$, i.e.\ the
length of the segment $XY$. The notation $XY$ will usually refer to the
line $XY$, and occasionally the ray or segment if explicitly stated.
\end{define}

\subsection*{Outline of Proof of Theorem~\ref{ismalldent}}

Theorem~\ref{ismalldent} (Theorem~\ref{smalldent}) states that the dented polygon described there
is a DDNA-polygon if and only if the angle $\alpha$ is small. For example,
the isosceles triangle in Figure~\ref{figPdelta} is a DDNA-polygon with
respect to its hypotenuse, but it is not a DDNA-polygon with respect to
its other sides.

Here is an overview of the proof of Theorem~\ref{smalldent}. A fundamental
idea is that if $\Gamma$ is a polygon that does not satisfy the DNA
Inequality, then for any polygonal DNA $\gamma$ providing a counterexample,
we can ``simplify'' it (using Lemmas \ref{lemma1}--\ref{wmacritical})
to obtain a special counterexample. For DNA of this special type,
the verification that the DNA Inequality fails or holds is  much easier.

Lemmas \ref{lemma1} and \ref{wmalemma1} clarify when we can assume that
the vertices of a counterexample are on an edge of $\Gamma$ or even
coincide with a vertex of $\Gamma$. If we could assume that all of the
vertices of any $CX_\Gamma$-polygon are vertices of $\Gamma$, then the
proof would be relatively easy. Since we cannot make this assumption,
we identify (Definition~\ref{separable}) manageable cells $\Gamma$,
which we call separable polygons, and a finite number of useful points
on the boundary of $\Gamma$ that are not vertices.  These points plus
the vertices of $\Gamma$ form the finite set $C$ of critical points.

The technical Lemma~\ref{wmacritical} shows that we may assume that the
vertices of our $CX_\Gamma$-polygon are all in the finite set $C$ or
else our $CX_\Gamma$-polygon has a special form involving vertices
of $\Gamma$ plus one or two other points on the boundary of $\Gamma$.
Thus Lemma~\ref{wmacritical} identifies special types of counterexamples
that any separable polygon $\Gamma$, that does not satisfy the DNA
Inequality, must contain.  Similar to the proof of Theorem~5.1 from
\cite{L-R}, we do this by removing ``jumps,'' i.e.\ adjacent pairs of
vertices such that the line segment connecting them is not contained in
the boundary of $\Gamma$. Such a line segment can intersect the boundary
just at its endpoints, in which case it is termed a ``jump''. (Since we
are trying to reduce to polygons having only critical vertices, we also
require that at least one of the endpoints is not a critical point for
it to be considered a jump.)  However, since the cell can be non-convex,
the interior of the line segment can intersect the boundary, in which
case we call it a ``leap''. The proof of Lemma~\ref{wmacritical} then
proceeds by descent: given a $CX_\Gamma$-polygon $\gamma$, we construct
a $CX_\Gamma$-polygon $\gamma'$ which either has a smaller jump number
(this is the number of jumps, with jumps such that neither endpoint
is a critical point double-counted), or the same jump number and fewer
leaps. To complete the proof of Lemma~\ref{wmacritical}, we again employ
the method of descent: we show that given a $CX_\Gamma$-polygon $\gamma$
with no jumps we can construct another $CX_\Gamma$-polygon $\gamma'$
having no jumps and fewer vertices which are not critical points.

Then, in Section~\ref{secseq} we turn to the problem of determining if
the DNA Inequality holds for our special types of polygons. If the DNA
Inequality does not hold for arbitrarily small dents of $P$, then we have
a sequence $P_{\delta_k}$ of dented polygons, and a sequence of counterexamples
$\gamma^k$, which we may assume to be of our special type.  The second
fundamental idea is to observe that all of these counterexamples have a
bounded number of vertices, and the set of all polygons contained in $P$
with a bounded number of vertices is compact; therefore, our sequence
of counterexamples has a limiting point, say $\overline{\gamma}$.
By studying which sequences can approach an equality case of the DNA
Inequality (in Lemma~\ref{seq}), we are able to show that infinitely
many of the $P_{\delta_k}$ contain counterexamples of a very special
type (which are described in Definition~\ref{gammakv}).  Potential
counterexamples of this type are so special that it is easy to specify
an inequality that determines whether they are indeed counterexamples,
i.e., whether $f_\Gamma < 0$.  From these inequalities, we employ
an analytic argument (in Section~\ref{secsmalldent}) to deduce
Theorem~\ref{smalldent}.

\section{\label{seclems} Three Useful Lemmas}

In this section, we give some useful machinery that will apply in any
cell $\Gamma$.

\begin{define} \label{freetomove}
We make the following important definitions concerning
modifications to the DNA polygon. 
\begin{itemize}
\item
If replacing $\gamma_i$ with any other point on the line
$\gamma_{i-1}\gamma_i$ sufficiently close to $\gamma_i$ yields a curve
contained within $\Gamma$, we say that $\gamma_i$ is \emph{free to move} along
the line $\gamma_{i-1}\gamma_i$.
\item
If $\gamma$ is a closed curve such that,
for all $i$, $\gamma_i$ is not free to move along $\gamma_{i-1}\gamma_i$
or $\gamma_i\gamma_{i+1}$, we say that $\gamma$ is a \emph{1-curve}.
\end{itemize}
\end{define}

\begin{lm} \label{lemma1} 
If $\gamma$ is a closed curve where there exists $i$ such that $\gamma_i$
is free to move along line $\gamma_{i-1}\gamma_i$, then one can always
move $\gamma_i$ one direction along $\gamma_{i-1}\gamma_i$, decreasing
$f_\Gamma$, until $\gamma_i$ becomes collinear with  $\gamma_{i+1},
\gamma_{i+2}$ in that order, or is no longer free to move.  If $\gamma_i$
is no longer free to move, then one of the following occurs:
\begin{itemize}
\item $\gamma_i$ reaches a vertex of $\Gamma$;
\item $\gamma_i$ reaches an edge of $\Gamma$ such that $\gamma_{i-1}$
does not lie on the line containing that edge;
\item the line segment $\gamma_i\gamma_{i+1}$ intersects the boundary at
a point other than $\gamma_{i+1}$ or $\gamma_i$.
\end{itemize}
\end{lm}

\begin{proof}

First note that the bulleted items simply give a list of possibilities,
such that it is necessary for one of them to hold if a vertex is no
longer free to move.  (Not all of them are sufficient.)

\smallskip

\noindent
We distinguish 2 cases:

\paragraph{Case 1:} 
$\gamma_{i-1}$ and $\gamma_{i+2}$ are on the same side of line
$\gamma_i\gamma_{i+1}$, as pictured in Figure~\ref{figcase1}.

\begin{figure}[ht]
\begin{center}
\unitlength 1.0 mm
\begin{picture}(32,20)(0,0)
\path(0,0)(10,10)(20,10)(30,0)
\put(-2.5,2){\makebox(0,0)[cc]{$\gamma_{i-1}$}}
\put(10,12){\makebox(0,0)[cc]{$\gamma_i$}}
\put(20,12){\makebox(0,0)[cc]{$\gamma_{i+1}$}}
\put(32,2){\makebox(0,0)[cc]{$\gamma_{i+2}$}}
\end{picture}
\end{center}
\caption{Diagram for Case 1} \label{figcase1}
\end{figure}
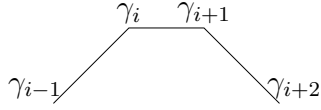

Moving $\gamma_i$ along line $\gamma_{i-1}\gamma_i$ in the direction
that increases the distance to $\gamma_{i-1}$ increases the perimeter,
but fixes the curvature, therefore decreasing $f_\Gamma$.

\paragraph{Case 2:}
They are on different sides, as pictured in Figure~\ref{figcase2}.

Let $H$ be the foot of the perpendicular from $\gamma_{i+1}$ to
line $\gamma_{i-1}\gamma_i$. Define $\theta$ to be angle 
$\angle H\gamma_{i+1}\gamma_i$. Let $a$ be the length of $H\gamma_{i+1}$.

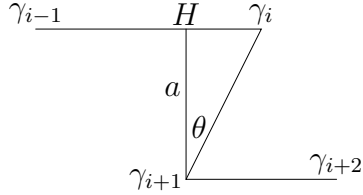
\begin{figure}[ht]
\begin{center}
\unitlength 1.0 mm
\begin{picture}(40,30)(0,0)
\path(0,20)(30,20)(20,0)(40,0)
\path(20,0)(20,20)
\put(21.75,7){\makebox(0,0)[cc]{$\theta$}}
\put(0,22){\makebox(0,0)[cc]{$\gamma_{i-1}$}}
\put(30,22){\makebox(0,0)[cc]{$\gamma_i$}}
\put(16,0){\makebox(0,0)[cc]{$\gamma_{i+1}$}}
\put(40,2){\makebox(0,0)[cc]{$\gamma_{i+2}$}}
\put(20,22){\makebox(0,0)[cc]{$H$}}
\put(18.25,12){\makebox(0,0)[cc]{$a$}}
\end{picture}
\end{center}
\caption{Diagram for Case 2} \label{figcase2}
\end{figure}

We will prove that $df_\Gamma/d\theta$ has at most one root
for $\theta \in (-\pi/2, \pi/2)$. 
\begin{align*}
0 = f_\Gamma' &= \alpha \cdot \frac{d}{d\theta}(\mbox{curvature}) - \frac{d}{d\theta} (\mbox{perimeter}) \\
&= 2\alpha - \frac{d}{d\theta}(a(\sec(\theta) + \tan(\theta))) \\
&= 2\alpha - \frac{a(1 + \sin(\theta))}{\cos^2(\theta)} \\
\Longleftrightarrow \frac{1 + \sin(\theta)}{\cos^2(\theta)} &= \frac{2\alpha}{a}
\end{align*}

Therefore, it suffices to show that 
$\frac{d}{d\theta}(\frac{1 + \sin(\theta)}{\cos^2(\theta)}) \ne 0$
on $(-\pi/2, \pi/2)$.
$$\frac{d}{d\theta}\left(\frac{1 + \sin(\theta)}{\cos^2(\theta)}\right) = \frac{(1 + \sin(\theta))^2}{\cos^3(\theta)} > 0 \mbox{ on } (-\pi/2, \pi/2).$$

Now, I claim that this finishes the proof of this Lemma. To see this,
observe that as $\theta \to \pi/2$, we have $f_\Gamma \to -\infty$.
Thus, as $f_\Gamma'$ has at most one root on $(-\pi/2, \pi/2)$, we either
have that $f_\Gamma$ is always decreasing, in which case we can move
$\gamma_i$ to the right, or that there exists $\beta \in (-\pi/2, \pi/2)$
such that $f_\Gamma$ is decreasing on $(\beta, \pi/2)$, and increasing
on $(-\pi/2, \beta)$.  In the latter case, we can move $\gamma_i$ to
the right if $\theta > \beta$ and to the left if $\theta < \beta$.

\end{proof}

\begin{lm} \label{wmalemma1}
If $\gamma$ is a $CX_\Gamma$-polygon, then there is a $CX_\Gamma$-polygon
$\gamma'$ which is a 1-curve.
\end{lm}

\begin{proof}
Assume there is some $CX_\Gamma$-polygon $\gamma$. Consider
$S = \{\gamma' \in F_\Gamma | l(\gamma') \leq l(\gamma)\}$, where $F_\Gamma$
is the set of all curves contained within $\Gamma$, and $l(\gamma)$
is the number of vertices of $\gamma$. $S$ is a non-empty (it contains
$\gamma$) compact set, and $f_\Gamma$ is a lower semi-continuous
function, so there is some $\gamma' \in S$ with $f_\Gamma(\gamma')$
minimal. Now, as $f_\Gamma(\gamma') \leq f_\Gamma(\gamma) < 0$,
$\gamma'$ is a $CX_\Gamma$-polygon. But, if $\gamma'$ were not a
1-curve, then by Lemma~\ref{lemma1}, there would exist a $\gamma''$
with $f_\Gamma(\gamma'') < f_\Gamma(\gamma')$, and this $\gamma''$
would be in $S$ because the proof of Lemma~\ref{lemma1} does not add
any vertices, providing a contradiction.
\end{proof}

\begin{lm} \label{vertexbound}
Write $V_\gamma$ for the set of vertices of $\gamma$. If $f_\Gamma(\gamma) < 0$ 
for some closed curve $\gamma$ contained within $\Gamma$, then there
is a curve $\gamma'$ which satisfies $V_{\gamma'} \subseteq V_\gamma$,
$f_\Gamma(\gamma') < 0$, and has number of vertices less than or equal
to $|V_\gamma|^2 - |V_\gamma|$.
\end{lm}

\begin{proof}
It suffices to show that for any closed curve $\gamma$ with at
least $|V_\gamma|^2 -|V_\gamma| + 1$ vertices, we can construct a
curve $\gamma'$ which satisfies $V_{\gamma'} \subseteq V_\gamma$,
$f_\Gamma(\gamma') < 0$ and has fewer vertices than $\gamma$.

Such a curve $\gamma$ has at least $|V_\gamma|^2 -|V_\gamma| + 1$ edges,
counting multiplicity. But, the number of edges without multiplicity is at
most $|V_\gamma|^2 -|V_\gamma|$, if we view our edges as directed. So by
the pigeonhole principle there is some directed edge repeated by $\gamma$,
i.e. there exists $i$ and $j$ with $i < j$ such that
$\gamma_i = \gamma_j, \gamma_{i+1} = \gamma_{j+1}$.
Now, consider the two curves
$\gamma^0 = \gamma_0, \gamma_1, \ldots, \gamma_i, \gamma_{j+1}, \gamma_{j+2}, \ldots, \gamma_n = \gamma_0$,
and $\gamma^1 = \gamma_{j}, \gamma_{i+1}, \gamma_{i+2}, \ldots, \gamma_{j-1}, \gamma_j = \gamma_j$,
as pictured in Figure~\ref{figi1j4}.

\begin{figure}[ht]
\unitlength 0.5 mm
\begin{center}
\begin{picture}(75,75)(0,0)
\path(10,40)(40,70)(70,40)(40,10)(10,40)
\path(10,40)(70,40)
\put(84,40){\makebox(0,0)[cc]{$\gamma_1 = \gamma_4$}}
\put(-4,40){\makebox(0,0)[cc]{$\gamma_2 = \gamma_5$}}
\put(40,74){\makebox(0,0)[cc]{$\gamma_3$}}
\put(40,4){\makebox(0,0)[cc]{$\gamma_0 = \gamma_6$}}
\end{picture}
\end{center}
\caption{Picture for $i = 1, j = 4$}
\label{figi1j4}
\end{figure}
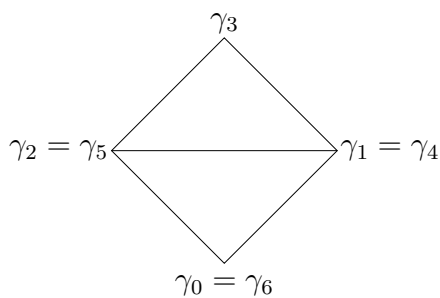

We have $f_\Gamma(\gamma^0) + f_\Gamma(\gamma^1) = f_\Gamma(\gamma) <0$,
so either $\gamma^0$ or $\gamma^1$ must satisfy the requirements
above for $\gamma'$.
\end{proof}

\section{\label{secsep} Separable Polygons}

This long section provides a detailed study of possible $CX_\Gamma$-polygons
in a special class of cells, which we term separable polygons.
We begin with the definitions.

\begin{define} \label{interior}
A vertex of a cell $\Gamma$ is called an \emph{interior vertex} if it is
contained in the interior of the convex hull of $\Gamma$.
\end{define}

\begin{define} \label{separable}
A polygon $\Gamma$ is called \emph{separable} if for any point $p$
in the interior of the cell determined by $\Gamma$, but not a vertex
of $\Gamma$, there is at most one interior vertex $v$ of $\Gamma$ such
that the straight line determined by $p$ and $v$ intersects $\Gamma$
in more than two points. Some examples of separable and non-separable
polygons are given in Figure~\ref{figsep}.

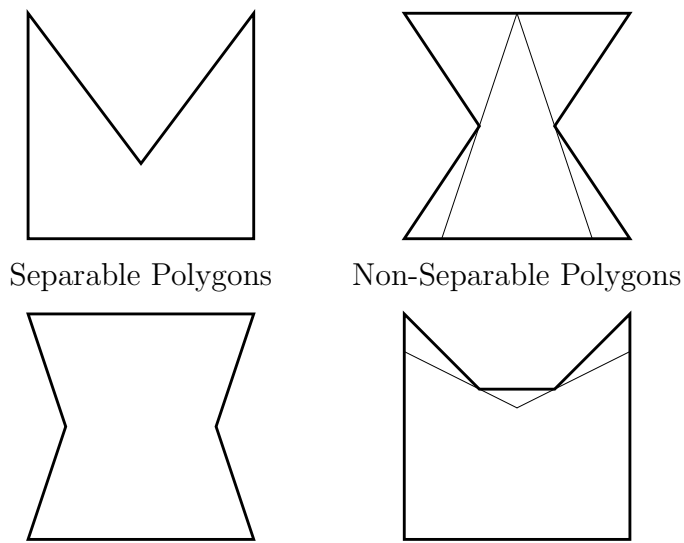
\begin{figure}[ht]
\begin{center}
\unitlength 1.0 mm
\allinethickness{0.3mm}
\begin{picture}(80,70)(0,0)
\Thicklines
\path(0,0)(5,15)(0,30)(30,30)(25,15)(30,0)(0,0)
\path(0,40)(0,70)(15,50)(30,70)(30,40)(0,40)
\path(50,0)(50,30)(60,20)(70,20)(80,30)(80,0)(50,0)
\path(50,40)(60,55)(50,70)(80,70)(70,55)(80,40)(50,40)
\thinlines
\path(65,70)(75,40)
\path(65,70)(55,40)
\path(65,17.5)(50,25)
\path(65,17.5)(80,25)
\put(15,35){\makebox(0,0)[cc]{Separable Polygons}}
\put(65,35){\makebox(0,0)[cc]{Non-Separable Polygons}}
\end{picture}
\end{center}
\caption{Examples of Separable and Non-Separable Polygons}
\label{figsep}
\end{figure}
\end{define}

\begin{cor2} \label{corinterior}
In a separable polygon $\Gamma$, if we assume that we have a
$CX_\Gamma$-polygon, then we have a $CX_\Gamma$-polygon, all of whose
vertices lie on the boundary.
\end{cor2}

\begin{proof}
Any vertex in the interior of a separable polygon is free to
move, along at least one of the two possible lines.  Therefore, by
Lemma~\ref{wmalemma1}, we may assume that we have a $CX_\Gamma$-polygon,
all of whose vertices lie on the boundary.
\end{proof}

\begin{define}
\label{critical}
The set of \emph{critical points} $C$ is the set of all vertices of $\Gamma$
plus any point $p$ in the interior of any edge of $\Gamma$ which is
collinear with two vertices of $\Gamma$, ($v, w$), which are distinct
from each other and the endpoints of the edge of $\Gamma$ upon which $p$
lies. Additionally, we require that the line segments connecting $pv,
pw$ are contained within $\Gamma$, and that $p$ is not free to move
along the line $pv$ (equivalently $pw$, as $p, v, w$ are collinear).
Figure~\ref{figexcrit} gives an example of a non-convex pentagon with $9$
critical points.

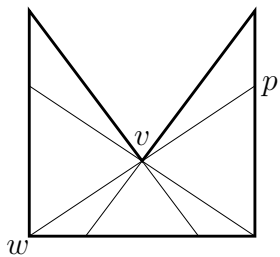
\begin{figure}[ht]
\begin{center}
\unitlength 1.0 mm
\begin{picture}(30,30)(0,0)
\Thicklines
\path(0,0)(0,30)(15,10)(30,30)(30,0)(0,0)
\thinlines
\path(30,30)(7.5,0)
\path(0,30)(22.5,0)
\path(0,0)(30,20)
\path(30,0)(0,20)
\put(32,20){\makebox(0,0)[cc]{$p$}}
\put(15,13){\makebox(0,0)[cc]{$v$}}
\put(-1.5,-1.5){\makebox(0,0)[cc]{$w$}}
\end{picture}
\end{center}
\caption{Examples of Critical Points}
\label{figexcrit}
\end{figure}
\end{define}

\begin{lm}\label{On2} $|C| \leq n^2$. \end{lm}
\begin{proof}
Consider a non-vertex critical point $p$ which is collinear with distinct
vertices $v, w$. Since $pv, pw$ are contained in $\Gamma$, so is $vw$. It
follows that $p$ must be the furthest point on the ray $vw$ such that
the line segment $vp$ is contained within $\Gamma$, or similarly for
ray $wp$. Thus, there are at most two non-vertex critical points for
each set $\{v, w\}$ of distinct vertices. So there are at most
$2 \cdot n(n-1) / 2$ non-vertex critical points, for a total of
$n(n-1) + n = n^2$ critical points in all.
\end{proof}

\begin{rem}
One can show that, if the cell $\Gamma$ is separable, $|C| \leq 2n -1$.
As we shall only need that it is bounded by a polynomial function
of $n$, we will leave the proof of this to an interested reader.
\end{rem}

\begin{define} \label{cutting}
If $\Gamma$ is separable, and $\gamma$ is a 1-curve, then we have
a way to split $\Gamma$ into two pieces, which we shall refer to as
\emph{cutting along segment $\gamma_{i-1}\gamma_i$.} We say that two
point $p$ and $q$, both in the interior of $\Gamma$ but not on line
segment $\gamma_{i-1}\gamma_i$, are on the same piece of $\Gamma$
if there is a (not necessarily closed) curve $c$ contained within
$\Gamma$ with endpoints $p$ and $q$ such that $c$ does not cross the
segment $\gamma_i\gamma_{i-1}$. We say that $c$ crosses the segment
$\gamma_i\gamma_{i-1}$ if there are two (possibly identical) points $x$
and $y$ which lie in the segment and on $c$ such that there are points
$x'$ and $y'$ on $c$, arbitrarily close to $x$ and $y$ respectively,
lying on opposite sides of the line $\gamma_i\gamma_{i+1}$. In
Figure~\ref{figcut}, two examples of this construction are given;
cutting along the bold line splits the region into the shaded parts and
the non-shaded parts.

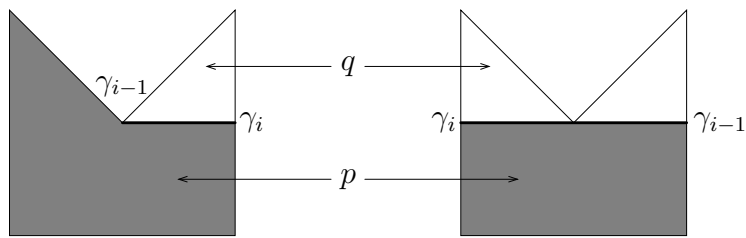
\begin{figure}[ht]
\begin{center}
\unitlength 1.5 mm
\begin{picture}(60,20)(0,0)
\thinlines
\path(40,10)(40,20)(50,10)(60,20)(60,10)
\shade\path(40,0)(40,10)(60,10)(60,0)(40,0)
\path(10,10)(20,20)(20,10)
\shade\path(0,0)(0,20)(10,10)(20,10)(20,0)(0,0)
\put(30,5){\makebox(0,0)[cc]{$p$}}
\put(30,15){\makebox(0,0)[cc]{$q$}}
\path(28.5,5)(15,5)
\path(28.5,15)(17.5,15)
\path(31.5,5)(45,5)
\path(31.5,15)(42.5,15)
\path(15.75,5.25)(15,5)(15.75,4.75)
\path(18.22,15.25)(17.5,15)(18.25,14.75)
\path(44.25,5.25)(45,5)(44.25,4.75)
\path(41.75,15.25)(42.5,15)(41.75,14.75)
\put(21.5,10){\makebox(0,0)[cc]{$\gamma_i$}}
\put(38.5,10){\makebox(0,0)[cc]{$\gamma_i$}}
\put(10,13.25){\makebox(0,0)[cc]{$\gamma_{i-1}$}}
\put(63,10){\makebox(0,0)[cc]{$\gamma_{i-1}$}}
\Thicklines
\path(40,10)(60,10)
\path(10,10)(20,10)
\end{picture}
\end{center}
\caption{Examples of Cutting Along $\gamma_{i-1}\gamma_i$}
\label{figcut}
\end{figure}

We also define a point $\gamma_i$ of the curve $\gamma$ to be a
``turn-around'' if $\gamma_{i+1}, \gamma_{i-2}$ lie on opposite
sides of $\gamma_{i-1}\gamma_i$, and when you cut along segment
$\gamma_{i-1}\gamma_i$, separating $\Gamma$ into two pieces,
$\gamma_{i+1}$ and $\gamma_{i-2}$ lie on different pieces (these are not
in general the same thing, as $\Gamma$ may be non-convex).  Additionally,
we require the same thing for $\gamma_{i-1}, \gamma_{i+2}$ with respect
to $\gamma_i\gamma_{i+1}$.  In Figure~\ref{figturn} is pictured a curve
$\gamma$ where $\gamma_i$ is a turn-around.

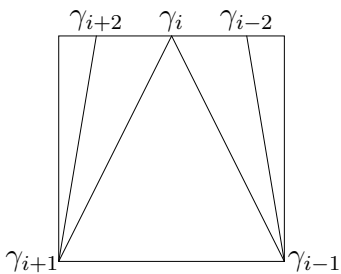
\begin{figure}[ht]
\begin{center}
\unitlength 1.0 mm
\begin{picture}(32.5,32.5)(0,0)
\path(0,0)(0,30)(30,30)(30,0)(0,0)
\path(0,0)(15,30)
\path(0,0)(5,30)
\path(30,0)(15,30)
\path(30,0)(25,30)
\put(15,32){\makebox(0,0)[cc]{$\gamma_i$}}
\put(34,0){\makebox(0,0)[cc]{$\gamma_{i-1}$}}
\put(-3.5,0){\makebox(0,0)[cc]{$\gamma_{i+1}$}}
\put(25,32){\makebox(0,0)[cc]{$\gamma_{i-2}$}}
\put(5,32){\makebox(0,0)[cc]{$\gamma_{i+2}$}}
\end{picture}
\end{center}
\caption{A ``Turn-Around''}
\label{figturn}
\end{figure}
\end{define}

\begin{lm}
\label{wmacritical}
In any separable polygon $\Gamma$ which contains a $CX_\Gamma$-polygon
$\gamma$, there exists a $CX_\Gamma$-polygon $\gamma'$ having one of
the following forms:
\begin{itemize}
\item $\gamma'$ has all vertices in $C$;
\item $\Gamma_0, \Gamma_1, \ldots, \Gamma_{i-1}, \Gamma_i, X, \Gamma_j, \Gamma_{j+1}, \ldots, \Gamma_n = \Gamma_0$,
where $i \leq j$;
\item $\Gamma_0, \Gamma_1, \ldots, \Gamma_{i-1}, \Gamma_i, X, Y, \Gamma_j, \Gamma_{j+1}, \ldots, \Gamma_n = \Gamma_0$,
where $i < j$;
\end{itemize}
where the vertices of $\Gamma$, in clockwise order, are
$\Gamma_0, \Gamma_1, \ldots, \Gamma_n = \Gamma_0$, $X$ is some point
on the boundary of cell $\Gamma$, and $Y$ is some point in the segment
$\Gamma_{j-1}\Gamma_j$ such that the line segment $XY$ intersects the
boundary in more than two points. In either case, we may assume that $X$
is a turn-around.
\end{lm}

\begin{proof}
We begin by assuming that there is some $CX_\Gamma$-polygon $\gamma$, but
there are no $CX_\Gamma$-polygons in the forms
\begin{gather*}
\Gamma_0, \Gamma_1, \ldots, \Gamma_{i-1}, \Gamma_i, X, \Gamma_j, \Gamma_{j+1}, \ldots, \Gamma_n = \Gamma_0 \mbox{\ or} \\
\Gamma_0, \Gamma_1, \ldots, \Gamma_{i-1}, \Gamma_i, X, Y, \Gamma_j, \Gamma_{j+1}, \ldots, \Gamma_n = \Gamma_0,
\end{gather*}
where $X$ is a turn-around, and prove that there is some
$CX_\Gamma$-polygon $\gamma'$ having vertices only in $C$.

In this proof, a ``jump'' is when we have $\gamma_i, \gamma_{i+1}$ which
are not both critical points such that the segment $\gamma_i\gamma_{i+1}$
intersects the boundary in {\it exactly} two points.  A jump is called a
``bad jump'' if neither of the $\gamma_i$ are critical points.  We term
the sum of the number of jumps and bad jumps (so bad jumps get counted
twice) the jump number of $\gamma$.  Additionally, we term a ``leap''
when we have $\gamma_i, \gamma_{i+1}$ such that $\gamma_i\gamma_{i+1}$
is not contained within the boundary of $\Gamma$.

\begin{cl} \label{lessleaps}
Suppose we have a $CX_\Gamma$-polygon (with at least one jump) that has
a leap $\gamma_i\gamma_{i+1}$ such that $\gamma_{i-1}, \gamma_{i+2}$
lie on the same side of line $\gamma_i\gamma_{i+1}$ or cutting along
$\gamma_i\gamma_{i+1}$ leaves $\gamma_{i-1}, \gamma_{i+2}$ on the
same piece of $\Gamma$. Then, there exists another $CX_\Gamma$-polygon
$\gamma'$ which either has a smaller jump number, or less leaps and the
same jump number.
\end{cl}

We first do the case when $\gamma_{i-1}, \gamma_{i+2}$ lie on the same
side of line $\gamma_i\gamma_{i+1}$.

Let $\gamma_i\gamma_{i+1}$ be a leap with $\gamma_{i-1}, \gamma_{i+2}$
on the same side of line $\gamma_i\gamma_{i+1}$.  We first examine
the case where $\gamma_i\gamma_{i+1}$ is not a jump, as pictured in
Figure~\ref{ii1nj}.  Observe that neither $\gamma_i$ nor $\gamma_{i+1}$
may be an interior vertex, as $\Gamma$ is separable.

\begin{figure}[ht]
\begin{center}
\unitlength 0.9 mm
\begin{picture}(105,50)(0,0)
\thinlines
\path(5,0)(15,40)(45,40)(35,20)(45,0)(5,0)
\path(90,5)(65,0)(60,0)(70,20)(60,40)(75,40)(85,30)
\put(8,3){\makebox(0,0)[cc]{$\gamma_i$}}
\put(10,37.5){\makebox(0,0)[cc]{$\gamma_{i+1}$}}
\put(68,3){\makebox(0,0)[cc]{$\gamma_i$}}
\put(70,37.5){\makebox(0,0)[cc]{$\gamma_{i+1}$}}
\Thicklines
\path(65,0)(75,40)(105,40)(95,20)(105,0)(65,0)
\path(30,5)(5,0)(0,0)(10,20)(0,40)(15,40)(25,30)
\put(20,20){\makebox(0,0)[cc]{$\gamma'$}}
\put(80,20){\makebox(0,0)[cc]{$\gamma^0$}}
\end{picture}
\end{center}
\caption{The case where $\gamma_i\gamma_{i+1}$ is not a jump} \label{ii1nj}
\end{figure}
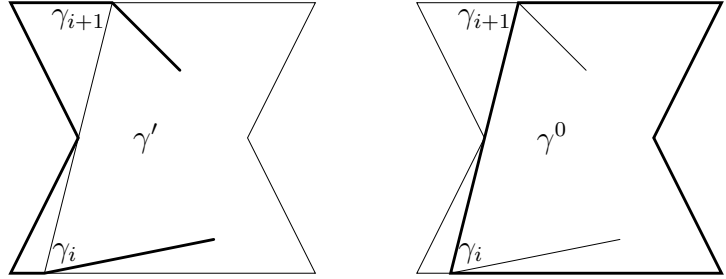

Define the closed curve $\gamma^0$ (pictured on the right) to be the
curve which consists of following the boundary of $\Gamma$, minus the
portion between $\gamma_i$ and $\gamma_{i+1}$, and jumping instead
from $\gamma_i \longrightarrow \gamma_{i+1}$.  Also, define closed curve
$\gamma'$ (a portion of which is pictured to the left) to be the curve
which consists of following $\gamma$, minus $\gamma_i \longrightarrow \gamma_{i+1}$,
and instead following the portion of the boundary which $\gamma^0$
does not follow.  Because of the orientation of the angles at
$\gamma_i, \gamma_{i+1}$ (which must be similar to as pictured above as
$\gamma_i, \gamma_{i+1}$ are not interior vertices and
$\gamma_{i-1}, \gamma_{i+2}$ lie on the same side of line
$\gamma_i\gamma_{i+1}$), we have
$f_\Gamma(\gamma) = f_\Gamma(\gamma^0) + f_\Gamma(\gamma')$.
(This equality uses $f_\Gamma(\Gamma) = 0$.)
Because $f_\Gamma(\gamma) < 0$, either $f_\Gamma(\gamma^0) < 0$ or
$f_\Gamma(\gamma') < 0$.  Now, both $\gamma', \gamma^0$ have a smaller
or equal jump number than $\gamma$ and fewer leaps (because our leap
was not a jump), so we may take one with a negative value of $f_\Gamma$
to be our $CX_\Gamma$-polygon $\gamma'$.

We now turn to the case where $\gamma_i\gamma_{i+1}$ is a jump such that
$\gamma_i$ is not free to move along line $\gamma_{i-1}\gamma_i$ and
$\gamma_{i+1}$ is not free to move along line $\gamma_{i+1}\gamma_{i+2}$,
as pictured in Figure~\ref{figlessleaps1}.

\begin{figure}[ht]
\begin{center}
\unitlength 0.75mm
\begin{picture}(145,50)(0,0)
\thinlines
\path(0,0)(0,30)(10,35)(15,45)(50,45)(55,35)(65,30)(65,0)(0,0)
\path(0,25)(20,45)
\path(65,25)(45,45)
\path(80,0)(80,30)(90,35)(95,45)(130,45)(135,35)(145,30)(145,0)(80,0)
\path(80,20)(145,20)
\put(57,37){\makebox(0,0)[cc]{$w$}}
\put(8,37){\makebox(0,0)[cc]{$v$}}
\put(17.5,37.5){\makebox(0,0)[cc]{$l_v$}}
\put(48,38){\makebox(0,0)[cc]{$l_w$}}
\put(-2.5,20){\makebox(0,0)[cc]{$\gamma_i$}}
\put(45,-2.5){\makebox(0,0)[cc]{$\gamma_{i+2}$}}
\put(20,-2.5){\makebox(0,0)[cc]{$\gamma_{i-1}$}}
\put(70,20){\makebox(0,0)[cc]{$\gamma_{i+1}$}}
\Thicklines
\path(100,0)(80,20)(80,25)(100,45)(125,45)(145,25)(145,20)(125,0)
\path(0,20)(65,20)
\path(0,20)(20,0)
\path(65,20)(45,0)
\end{picture}
\end{center}
\caption{First Figure where $\gamma_i\gamma_{i+1}$ is a jump}
\label{figlessleaps1}
\end{figure}
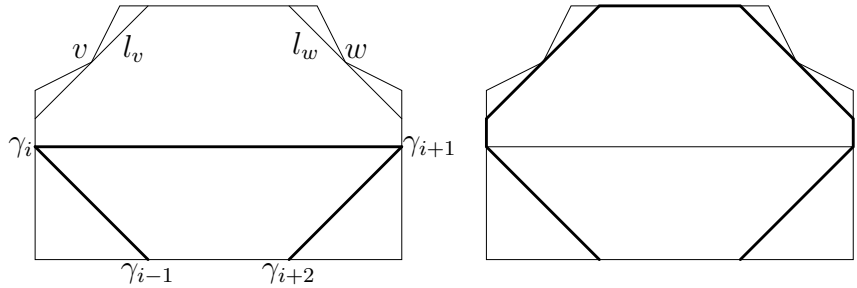

If we cut along the line $\gamma_i\gamma_{i+1}$, this separates $\Gamma$
into two pieces.  On the side not containing $\gamma_{i+2}, \gamma_{i-1}$,
through every interior vertex $v$ of $\Gamma$, we construct a line $l_v$
which passes through $v$ but does not intersect the interior of the line
segment $\gamma_i\gamma_{i+1}$.  Now, we form a new curve $\gamma'$
by replacing $\gamma_i\gamma_{i+1}$ with the path that goes along the
boundary and the $l_v$, as pictured above.  The resulting curve has the
same curvature (the orientation of the angles at $\gamma_i$ is similar
to as pictured above because $\gamma_i$ is not free to move along line
$\gamma_{i-1}\gamma_i$; similarly, the possible orientations of angles at
$\gamma_{i+1}$ are limited), but a greater perimeter, and thus a smaller
(and hence negative) value of $f_\Gamma$, while having one less jump,
completing the proof of this case.

Now, we consider the case when (without loss of generality)
$\gamma_i\gamma_{i+1}$ is a jump and $\gamma_i$ is free to move
along line $\gamma_{i-1}\gamma_i$. In this case, move it along
the line until it is no longer free to move; call the position
that it reaches $\gamma_i'$. Because $\Gamma$ is separable, line
$\gamma_i'\gamma_{i+1}$ must intersect the boundary at only two points
(the line $\gamma_{i-1}\gamma_i'$ intersects it in more than two). Thus,
$\gamma_i'$ must have reached the boundary. If $\gamma_{i+1}$ is free
to move along line $\gamma_{i+1}\gamma_{i+2}$, then we construct in
a similar manner $\gamma_{i+1}'$ (otherwise, define $\gamma_{i+1}' = \gamma_{i+1}$).
Now, the curve formed by using $\gamma_i', \gamma_{i+1}'$
instead of $\gamma_i, \gamma_{i+1}$ has the same number of jumps,
one of which is $\gamma_i'\gamma_{i+1}'$.  But, by the previous case,
we can create a new curve without that jump.  Thus, this completes the
proof of the case when $\gamma_{i-1}, \gamma_{i+2}$ lie on the same side
of line $\gamma_i\gamma_{i+1}$.

Next, we do the case when $\gamma_{i-1}, \gamma_{i+2}$ lie on
different sides of $\gamma_i\gamma_{i+1}$, but cutting along
$\gamma_i\gamma_{i+1}$ leaves them on the same piece of $\Gamma$, as
pictured in Figure~\ref{figlessleaps2}.

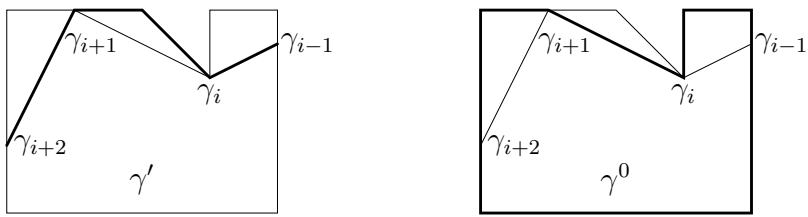
\begin{figure}[ht]
\begin{center}
\unitlength 0.9 mm
\begin{picture}(110,30)(0,0)
\thinlines
\path(0,0)(0,30)(10,30)(30,20)(30,30)(40,30)(40,0)(0,0)
\path(70,10)(80,30)(90,30)(100,20)(110,25)
\put(20,5){\makebox(0,0)[cc]{$\gamma'$}}
\put(90,5){\makebox(0,0)[cc]{$\gamma^0$}}
\multiput(44.5,25)(70,0){2}{\makebox(0,0)[cc]{$\gamma_{i-1}$}}
\multiput(30,17.5)(70,0){2}{\makebox(0,0)[cc]{$\gamma_i$}}
\multiput(12.5,25)(70,0){2}{\makebox(0,0)[cc]{$\gamma_{i+1}$}}
\multiput(5,10)(70,0){2}{\makebox(0,0)[cc]{$\gamma_{i+2}$}}
\Thicklines
\path(0,10)(10,30)(20,30)(30,20)(40,25)
\path(70,0)(70,30)(80,30)(100,20)(100,30)(110,30)(110,0)(70,0)
\end{picture}
\end{center}
\caption{Second Figure where $\gamma_i\gamma_{i+1}$ is a jump}
\label{figlessleaps2}
\end{figure}

Define the closed curve $\gamma^0$ (pictured on the right) to be the curve
which consists of following the boundary of $\Gamma$, minus the portion
between $\gamma_i$ and $\gamma_{i+1}$, and jumping instead from
$\gamma_i \longrightarrow \gamma_{i+1}$, and define closed curve $\gamma'$ (a
portion of which is pictured to the left) to be the curve which consists
of following $\gamma$, minus $\gamma_i \longrightarrow \gamma_{i+1}$, and
instead following the portion of the boundary which $\gamma^0$ does
not follow.  Because of the orientation of the angles at
$\gamma_i, \gamma_{i+1}$ (which must be similar to as pictured above), we
have $f_\Gamma(\gamma) = f_\Gamma(\gamma^0) + f_\Gamma(\gamma')$.
Because $f_\Gamma(\gamma) < 0$, either $f_\Gamma(\gamma^0) < 0$ or
$f_\Gamma(\gamma') < 0$.  If $f_\Gamma(\gamma') < 0$, then we are done.
Otherwise, as the point of $\gamma^0$ which coincides with $\gamma_{i+1}$
is free to move, we can move it until it coincides with a critical point,
forming a curve with a negative value of $f_\Gamma$ with no jumps,
completing the proof of this claim.

\begin{cl} \label{jump0}
We may construct a curve $\gamma'$ with a jump number of 0,
which is also a $CX_\Gamma$-polygon.
\end{cl}

It clearly suffices to show that given a $CX_\Gamma$-polygon $\gamma$,
we can construct another $CX_\Gamma$-polygon with a smaller jump number,
or with the same jump number but having fewer leaps. Consider some jump
$\gamma_i\gamma_{i+1}$. Without loss of generality, let $\gamma_i$ not be
a critical point. By applying claim~\ref{lessleaps}, we may assume that for any $j$ such
that $\gamma_j\gamma_{j+1}$ is a leap, $\gamma_{j-1}$ and $\gamma_{j+2}$
do not lie on the same side of line $\gamma_j\gamma_{j+1}$. We may
further assume that for any $j$ such that $\gamma_j\gamma_{j+1}$ is a
leap, when we cut along segment $\gamma_j\gamma_{j+1}$, $\gamma_{j-1}$
and $\gamma_{j+2}$ are on different pieces.

If $\gamma_{i-1}\gamma_i$ is not a leap, then $\gamma_{i-1}$ must be on
the same edge of $\Gamma$ as $\gamma_i$. By Lemma~\ref{lemma1}, we can
move $\gamma_i$ along the line $\gamma_{i-1}\gamma_i$, until one of the
following occurs:
\begin{itemize}
\item It reaches a vertex: In this case, $\gamma_i\gamma_{i+1}$ either
is no longer a jump if $\gamma_{i-1}$ is a critical point or no longer
a bad jump otherwise; either way, the jump number decreases.
\item The line segment $\gamma_i\gamma_{i+1}$ intersects the boundary
at a point other than $\gamma_{i+1}$: Thus $\gamma_i\gamma_{i+1}$ is no
longer a jump, decreasing the jump number.
\item $\gamma_i$ becomes collinear with  $\gamma_{i+1}, \gamma_{i+2}$
in that order: This implies that we drop $\gamma_{i+1}$, and
$\gamma_i\gamma_{i+2}$ is not a jump, because it intersects the boundary
of $\Gamma$ in a third point (the previous location of $\gamma_{i+1}$).
\end{itemize}
Thus, we may assume that $\gamma_{i-1}\gamma_i$ is also a leap.  From our
earlier discussion about leaps, we may assume that $\gamma_{i+2}$ and
$\gamma_{i-1}$ lie on opposite sides of line $\gamma_i\gamma_{i+1}$,
and that $\gamma_{i-2}$ and $\gamma_{i+1}$ lie on opposite sides of line
$\gamma_{i-1}\gamma_i$, which is pictured in Figure~\ref{figjump0}.

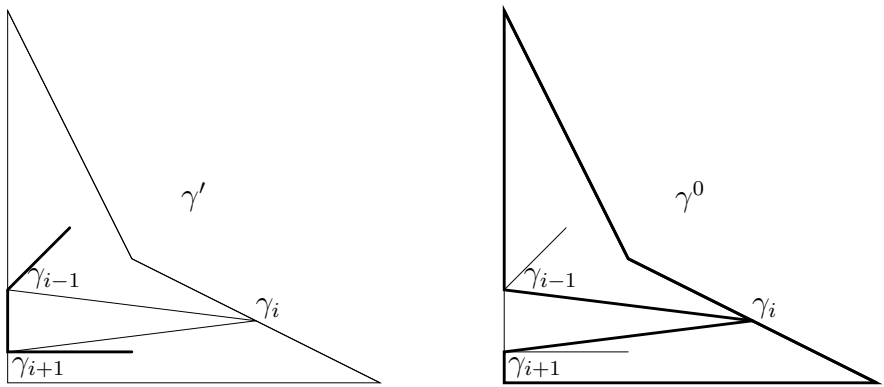
\begin{figure}[ht]
\begin{center}
\unitlength 1.65mm
\begin{picture}(70,30)(0,0)
\put(15,15){\makebox(0,0)[cc]{$\gamma'$}}
\put(55,15){\makebox(0,0)[cc]{$\gamma^0$}}
\thinlines
\path(20,5)(0,2.5)(0,0)(30,0)(10,10)(0,30)
\path(30,0)(10,10)(0,30)(0,7.5)(20,5)
\path(50,2.5)(40,2.5)(40,7.5)(45,12.5)
\Thicklines
\path(10,2.5)(0,2.5)(0,7.5)(5,12.5)
\path(60,5)(40,2.5)(40,0)(70,0)(50,10)(40,30)
\path(70,0)(50,10)(40,30)(40,7.5)(60,5)
\put(61,6){\makebox(0,0)[cc]{$\gamma_i$}}
\put(21,6){\makebox(0,0)[cc]{$\gamma_i$}}
\put(3.75,8.4){\makebox(0,0)[cc]{$\gamma_{i-1}$}}
\put(43.75,8.4){\makebox(0,0)[cc]{$\gamma_{i-1}$}}
\put(2.5,1.25){\makebox(0,0)[cc]{$\gamma_{i+1}$}}
\put(42.5,1.25){\makebox(0,0)[cc]{$\gamma_{i+1}$}}
\end{picture}
\end{center}
\caption{Figure for Claim~\ref{jump0}}
\label{figjump0}
\end{figure}

In this case, define the closed curve $\gamma^0$ (pictured on the right)
to be the curve which consists of following the boundary of $\Gamma$,
minus the portion between $\gamma_{i-1}$ and $\gamma_{i+1}$ (the portion
not containing $\gamma_i$), and jumping instead from
$\gamma_{i-1} \longrightarrow \gamma_i \longrightarrow \gamma_{i+1}$, and
define closed curve $\gamma'$ (a portion of which is pictured to the
left) to be the curve which consists of following $\gamma$, minus
$\gamma_{i-1} \longrightarrow \gamma_i \longrightarrow \gamma_{i+1}$, and
instead following the portion of the boundary which $\gamma^0$ does
not follow.  As $\gamma_{i \pm 2}, \gamma_{i \mp 1}$ lie on opposite
sides of $\gamma_i\gamma_{i \pm 1}$ and end up on different pieces when
we cut along segments $\gamma_i\gamma_{i \pm 1}$, the angles must be
oriented in a similar fashion to the ones in the above diagram, and we
thus have $f_\Gamma(\gamma) = f_\Gamma(\gamma^0) + f_\Gamma(\gamma')$.
Now, consider moving the vertices of $\gamma^0$ that coincide with
$\gamma_{i \pm 1}$ (not along the line connecting them to $\gamma_i$, but
along the other of two possible lines), until they become collinear with
$\gamma_i$, or until each one reaches vertices of $\Gamma$ or the line
segment joining that point to $\gamma_i$ is not a jump, constructing
a curve ${\gamma^0}'$. By Lemma~\ref{lemma1},
$f_\Gamma(\gamma^0) \geq f_\Gamma({\gamma^0}')$. Now, if
$f_\Gamma({\gamma^0}') \geq 0$, this implies $f_\Gamma(\gamma') < 0$.
But, $\gamma'$ has two less jumps than $\gamma$.  On the other hand if
$f_\Gamma({\gamma^0}') < 0$, then if the point coinciding with
$\gamma_{i \pm 1}$ became collinear with $\gamma_i$, the curve
${\gamma^0}'$ has no jumps and has a negative value of $f_\Gamma$,
completing the proof of this claim.  Otherwise, by assumption (see first
paragraph of the proof), one of the following holds:
\begin{itemize}
\item $\gamma_i$ is not a turn-around in ${\gamma^0}'$: Therefore, we
may apply claim~\ref{lessleaps} to the curve ${\gamma^0}'$ to produce a curve which
either has a smaller jump number or an identical jump number but fewer
leaps than ${\gamma^0}'$.
\item If we write
${\gamma^0}' = \Gamma_0, \Gamma_1, \ldots, \Gamma_{i-1}, \Gamma_i, X, \Gamma_j, \Gamma_{j+1}, \ldots, \Gamma_n = \Gamma_0$,
we have $i > j$: It is clear that
$f_\Gamma({\gamma^0}') = f_\Gamma(X, \Gamma_j, \Gamma_{j+1}, \ldots, \Gamma_i, X)$.
But to the latter curve, we may apply claim~\ref{lessleaps}, to produce a curve which
either has a smaller jump number or an identical jump number but fewer
leaps than ${\gamma^0}'$.
\item If we write
${\gamma^0}' = \Gamma_0, \Gamma_1, \ldots, \Gamma_{i-1}, \Gamma_i, X, Y, \Gamma_j, \Gamma_{j+1}, \ldots, \Gamma_n = \Gamma_0$,
we have $i \geq j$: In this case, exactly the same argument works, replacing
$X, \Gamma_j, \Gamma_{j+1}, \ldots, \Gamma_i, X$ with
$X, Y, \Gamma_j, \Gamma_{j+1}, \ldots, \Gamma_i, X$.
\end{itemize}
This completes the proof of this claim.

\smallskip

As noted at the beginning of the proof, the next claim will complete
the proof of Lemma~\ref{wmacritical}:

\begin{cl} \label{j0tolm}
Given a $CX_\Gamma$-polygon $\gamma$ with a jump number of 0, we may
construct another $CX_\Gamma$-polygon $\gamma'$ which consists of vertices
only in $C$.
\end{cl}

It clearly suffices to show that given such a $CX_\Gamma$-polygon
$\gamma$, we can construct another $CX_\Gamma$-polygon $\gamma'$ with
less vertices not in $C$, which also has a jump number of 0.  I first
claim that we may assume that no vertex not in $C$ is free to move.
For if any are free to move, then we may move them until that is no longer
the case, and we will not increase the number of vertices in $C$. It is
clear that this operation cannot increase the jump number.  Now, take some
$\gamma_i \notin C$. As it is not free to move, we have without loss of
generality, $\gamma_i, X, \gamma_{i+1}$ collinear in that order, for some
interior vertex $X$ of $\Gamma$. Of course, we cannot have $\gamma_{i+1}$
a critical point either, as that would imply $\gamma_i$ is a vertex.

\paragraph{Case 1:}
$\gamma_{i+2} \neq \gamma_i$ and $\gamma_{i+1} \neq \gamma_{i-1}$.
It follows that $\gamma_{i+2}, \gamma_{i+1}$ lie on the same edge of
$\Gamma$, as do $\gamma_i, \gamma_{i-1}$.  Figure~\ref{figorient}
shows the three ways that these pairs of vertices can lie on their
respective edges.

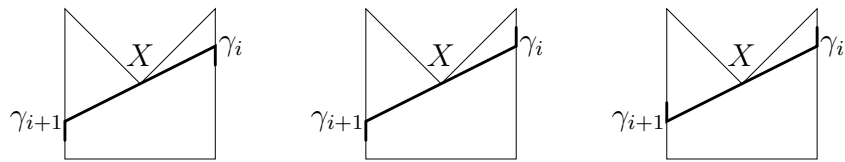
\begin{figure}[ht]
\begin{center}
\unitlength 1.0 mm
\begin{picture}(100,20)(0,0)
\thinlines
\path(0,0)(0,20)(10,10)(20,20)(20,0)(0,0)
\path(40,0)(40,20)(50,10)(60,20)(60,0)(40,0)
\path(80,0)(80,20)(90,10)(100,20)(100,0)(80,0)
\Thicklines
\path(0,2.5)(0,5)(20,15)(20,12.5)
\path(40,2.5)(40,5)(60,15)(60,17.5)
\path(80,7.5)(80,5)(100,15)(100,17.5)
\multiput(10,13.8)(40,0){3}{\makebox(0,0)[cc]{$X$}}
\multiput(22.1,15)(40,0){3}{\makebox(0,0)[cc]{$\gamma_i$}}
\multiput(-3.75,5)(40,0){3}{\makebox(0,0)[cc]{$\gamma_{i+1}$}}
\end{picture}
\end{center}
\caption{Possible Orientations for $\gamma_{i- 1},\gamma_i,\gamma_{i+1},\gamma_{i+2}$}
\label{figorient}
\end{figure}

In this case, consider rolling the line $\gamma_i\gamma_{i+1}$ around
$X$, as pictured in Figure~\ref{figrolling}. As a function of the angle
$\theta$ that line $\gamma_i\gamma_{i+1}$ makes with some fixed line,
I next show that $f_\Gamma$ is concave down, at least for the angles for
which $\gamma_i, \gamma_{i+1}$ remain on the same edge of $\Gamma$ and on
the same side of line $\gamma_{i-1}\gamma_{i+2}$ as they were originally.

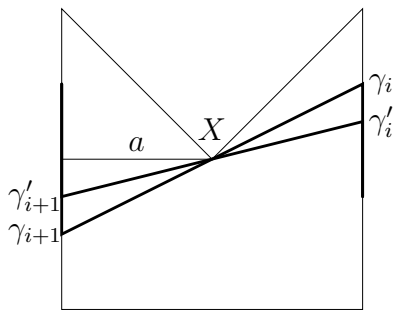
\begin{figure}[ht]
\begin{center}
\unitlength 1.0 mm
\begin{picture}(40,40)(0,0)
\thinlines
\path(0,0)(0,40)(20,20)(40,40)(40,0)(0,0)
\path(0,20)(20,20)
\Thicklines
\path(0,30)(0,10)(40,30)(40,15)
\path(0,30)(0,15)(40,25)(40,15)
\put(42.5,30){\makebox(0,0)[cc]{$\gamma_i$}}
\put(42.5,25){\makebox(0,0)[cc]{$\gamma_i'$}}
\put(-3.5,10){\makebox(0,0)[cc]{$\gamma_{i+1}$}}
\put(-3.5,15){\makebox(0,0)[cc]{$\gamma_{i+1}'$}}
\put(20,24){\makebox(0,0)[cc]{$X$}}
\put(10,22){\makebox(0,0)[cc]{$a$}}
\end{picture}
\end{center}
\caption{``Rolling'' about $X$}
\label{figrolling}
\end{figure}

Recall that $f_\Gamma(\gamma) = \alpha \cdot (\mbox{curvature}) -
(\mbox{perimeter})$. For $\theta$ in the interval specified above, the
curvature is clearly linear, so it suffices to show that the perimeter
function is concave up. Now, the perimeter of $\gamma$ is a constant plus
the sum of lengths
$\dist(\gamma_{i-1},\gamma_i) + \dist(\gamma_i,X) + \dist(X,\gamma_{i+1}) + \dist(\gamma_{i+1},\gamma_{i+2})$.
Thus, by symmetry, it suffices to show that
$\dist(\gamma_{i-1},\gamma_i) + \dist(\gamma_i,X)$ is a concave-up
function of $\theta$.  This clearly does not depend on the choice
of our fixed line, so we let our fixed line be the perpendicular
from $X$ to the edge of $\Gamma$ upon which $\gamma_i, \gamma_{i-1}$
lie. Then, for $\theta$ in the above domain, depending upon orientation,
$\dist(\gamma_{i-1},\gamma_i) + \dist(\gamma_i,X)$ is given up to a
constant by: $a(\sec{\theta} \pm \tan{\theta})$, where $a$ is the length
of the perpendicular from $X$ to that side. The second derivative of that
expression is given by $\frac{\cos{\theta}}{(1 \mp \sin{\theta})^2} > 0$
for $\theta$ in that domain, since that domain is always contained in
$(-\pi/2, \pi/2)$. Thus, $f_\Gamma$ is a concave-down function in that
domain, so the minimum of $f_\Gamma$ as we roll our line around $X$ occurs
at the end points of the domain. If we replace $\gamma$ by the curve that
uses this minimum instead, we have not increased the number of vertices
which are not in $C$, nor have we increased the jump number, and we have
decreased the number of $\gamma_i$ which fall under this case.  Thus,
if there is some $\gamma_i \notin C$, we may assume that
$\gamma_{i+2} = \gamma_i$ or $\gamma_{i+1} = \gamma_{i-1}$.

\paragraph{Case 2:}
$\gamma_{i+1} = \gamma_{i-1}$, but we have $\gamma_i \neq$ both $\gamma_{i \pm 2}$,
from which it follows that $\gamma_{i + 2}$ lies on the same
edge as $\gamma_{i + 1}$, and $\gamma_{i - 2}$ lies on the same edge as
$\gamma_{i - 1}$, as pictured in Figure~\ref{figj0tolmc2}.

\begin{figure}[ht]
\begin{center}
\unitlength 1.0 mm
\begin{picture}(100,20)(0,0)
\thinlines
\path(0,0)(0,20)(10,10)(20,20)(20,0)(0,0)
\path(40,0)(40,20)(50,10)(60,20)(60,0)(40,0)
\path(80,0)(80,20)(90,10)(100,20)(100,0)(80,0)
\Thicklines
\path(0,7.5)(20,12.5)(20,2.5)
\path(40,7.5)(60,12.5)(60,17.5)
\path(80,7.5)(100,12.5)(100,2.5)
\path(80,7.5)(100,12.5)(100,17.5)
\multiput(10,13.8)(40,0){3}{\makebox(0,0)[cc]{$X$}}
\multiput(-2,7.5)(40,0){3}{\makebox(0,0)[cc]{$\gamma_i$}}
\multiput(23.85,12.5)(40,0){3}{\makebox(0,0)[cc]{$\gamma_{i \pm 1}$}}
\end{picture}
\end{center}
\caption{Figure for Case 2 (of Claim~\ref{j0tolm})}
\label{figj0tolmc2}
\end{figure}
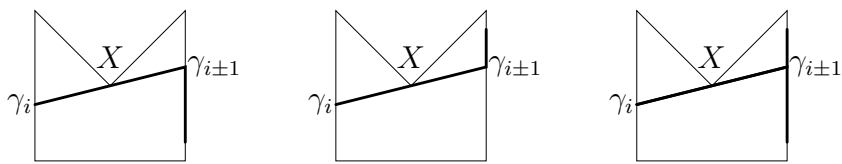

In this case, consider rolling the lines $\gamma_i\gamma_{i+1}$
and $\gamma_{i-1}\gamma_i$ around $X$ together, so that we keep
$\gamma_{i-1} = \gamma_{i+1}$. Similar to the previous case,
we will show that $f_\Gamma$ is concave down (in the appropriate
interval). Again, the curvature is linear, so it suffices to show
that the perimeter is concave up. The perimeter, up to an additive
constant, is given by
$2\dist(\gamma_i,X) + \dist(X,\gamma_{i+1}) + \dist(\gamma_{i+1},\gamma_{i+2}) + \dist(\gamma_{i-1},X) + \dist(\gamma_{i-2},\gamma_{i-1})$.
The calculation in Case 1 showed that
$\dist(X,\gamma_{i+1}) + \dist(\gamma_{i+1},\gamma_{i+2})$ and
$\dist(\gamma_{i-1},X) + \dist(\gamma_{i-2},\gamma_{i-1})$ are concave up,
so it suffices to show that $\dist(\gamma_i,X)$ is a concave-up function
of $\theta$. Again, choosing our fixed line to be from $X$ to the edge of
$\Gamma$ upon which $\gamma_i$ lies, we see that our function is given,
up to a constant, by $a\sec{\theta}$, which is concave up in $(-\pi/2, \pi/2)$.
Thus, the minimum of $f_\Gamma$, as we roll our line around
$X$, occurs at the end points of the domain. As in the previous case,
we may assume that there is some $\gamma_i$ which does not fall under
this case or the previous one, provided that, after this reduction and
the previous one, we still have some $\gamma_i \notin C$.

\paragraph{Case 3:}
$\gamma_{i+1} = \gamma_{i-1}$, and $\gamma_i = \gamma_{i+2}$ or $\gamma_{i-2}$.
Without loss of generality, say that $\gamma_i = \gamma_{i+2}$.
Now, define the curves
\begin{align*}
\gamma^0 &:= \gamma_0, \gamma_1, \ldots, \gamma_{i-1}, \gamma_{i+2}, \ldots \gamma_n = \gamma_0 \\
\gamma^1 &:= \gamma_{i-1}, \gamma_i, \gamma_{i+1} = \gamma_{i-1}.\end{align*}
We have $f_\Gamma(\gamma^0) + f_\Gamma(\gamma^1) = f_\Gamma(\gamma) < 0$.
Now, consider rotating $\gamma^1$ around $X$ until one of its
vertices becomes equal to a vertex of $\Gamma$, producing a new
curve ${\gamma^1}'$. In order to prove that we can do this to decrease
$f_\Gamma$, it suffices to show that, in terms of the angle, $f_\Gamma$ is
concave down. As the curvature is constant, it suffices to show that the
perimeter is concave up. As the perimeter is given by
$2(\dist(\gamma_i,X) + \dist(\gamma_{i+1},X))$, we have already
seen in Case 2 that this is concave up. Thus, we can construct a
curve ${\gamma^1}'$ that has all vertices in $C$, and a $\gamma^0$
that has fewer vertices not in $C$ than $\gamma$, such that
$f_\Gamma(\gamma^0) + f_\Gamma(\gamma^1) < 0, f_\Gamma({\gamma^1}') \leq f_\Gamma(\gamma^1)$.
This gives
$f_\Gamma({\gamma^1}') + f_\Gamma(\gamma^0) \leq f_\Gamma(\gamma^1) + f_\Gamma(\gamma^0) < 0$.
Thus, either $f_\Gamma({\gamma^1}') < 0$ or $f_\Gamma(\gamma^0) < 0$;
either way, we have constructed another curve with fewer vertices
$\notin C$ which is also a $CX_\Gamma$-polygon.

\smallskip

This completes the proof of this case, hence of this claim, and hence of this Lemma.
\end{proof}

\begin{cor} \label{lenbound}
If a cell $\Gamma$ is a separable polygon, then $\Gamma$ satisfies the
DNA Inequality if and only if the inequality holds when the DNA has
$|C|^2 - |C|$ or fewer vertices.
\end{cor}

\begin{proof}
Apply Lemmas \ref{vertexbound} and \ref{wmacritical};
note that $n + 2 \leq |C|^2 - |C|$. \end{proof}

\section{\label{secpoly} Proof of Theorem~\ref{ipoly}}

Here, we provide some estimates for the complexity of determining if a
given polygon satisfies the DNA Inequality.
For the proof of this theorem, we adopt a simplified model of computation,
described below in Definition~\ref{elem}.

\begin{define} \label{elem}
The following operations count as \emph{elementary operations}:
\begin{itemize}
\item Computing the sum, product, difference, or quotient of any two numbers.
\item Computing any trigonometric functions of any number.
\item Comparing two numbers (i.e.\ testing which one is larger).
\item Computing the roots of a polynomial equation, given its coefficients
(if the degree of the polynomial is bounded).
\end{itemize}
\end{define}

\begin{thm} \label{poly} 
There exists an algorithm which when given as input a separable 
polygon $\Gamma$ specified by its $n$ vertices,
determines whether or not $\Gamma$ satisfies the DNA 
inequality, in number of elementary operations polynomial in $n$.
\end{thm}

\begin{proof}
By Lemma~\ref{wmacritical}, to determine whether $\Gamma$ satisfies the DNA
Inequality, it suffices to examine curves with $|C|^2 - |C|$ or fewer
vertices, all in $C$, as well as the curves
\begin{equation}\label{specialcase}
\begin{split}
\Gamma_0, \Gamma_1, \ldots, \Gamma_{i-1}, \Gamma_i, X, \Gamma_j, \Gamma_{j+1}, \ldots, \Gamma_n = \Gamma_0 \mbox{\ and} \\
\Gamma_0, \Gamma_1, \ldots, \Gamma_{i-1}, \Gamma_i, X, Y, \Gamma_j, \Gamma_{j+1}, \ldots, \Gamma_n =  \Gamma_0.
\end{split}
\end{equation}
For these curves, there are a finite number of ways to choose $i, j$,
the edges upon which $X$ (and $Y$ if we are in the latter case) lie, and,
if we are in the latter case, the interior vertex which line $XY$ passes
through. For each combination, $f_\Gamma$ as a function of the position
of $X$ can be differentiated, and the curves can be considered for each
zero of the derivative and at the end points. For the curves in the form
(\ref{specialcase}) there are $O(n)$ ways to choose each of $i, j$, the
edges upon which $X$ (and possibly $Y$ as well) lie, and the interior
vertex on line $XY$ (if we are in the second case). As computing the
zeros of the derivative\footnote{To do this, we need only to compute
sums, differences, products, quotients, trigonometric functions, and
solve polynomials of bounded degree.} and the average curvature is linear
time in $n$, this part gives contribution $O(n^6)$ to the run time for
examining the second case, and $O(n^4)$ for examining the first case.

Thus, it suffices to show that one can check the curves with $|C|^2 -
|C|$ or fewer vertices, all in $C$ in time $O(n^{12}\log{n})$. By
Lemma~\ref{On2}, we have $|C| = O(n^2)$.  Define $S$ to the set of
ordered pairs of critical points such that the segment connecting them
lies within $\Gamma$.  Define the functions $f^k : S^2 \to \mathbb{R}$ of
$(e_1, e_2)$ to be the minimal possible value of the function $f_\Gamma$
over all (possibly open) polygonal paths with at most $k+2$ vertices,
whose first edge is $e_1$ and whose last edge is $e_2$ (if there are no
such polygonal paths, we assign value $\infty$).  Now, we can precompute
a table of values for $f^k$ for any $k$.

Suppose we want to find the minimal value that is assumed by all closed
curves with  $|C|^2 - |C|$ or fewer vertices.  If we assume that the
curve has three consecutive vertices $v_1, v_2, v_3$, then the minimal
value of $f_\Gamma$ for such a curve is
$$f^{|C|^2 - |C| - 1}(v_2v_1, v_3v_2) + \alpha(\pi - \angle{v_1v_2v_3}),$$
where $1/\alpha$ is the average curvature of $\Gamma$ (as in
Definition~\ref{f}). Thus, the DNA Inequality holds in $\Gamma$
if and only if
$$f^{|C|^2 - |C| - 1}(v_2v_1, v_3v_2) + \alpha(\pi - \angle{v_1v_2v_3}) > 0$$
for any $v_1, v_2, v_3$ such that $v_1v_2, v_2v_3 \in S$. So, if we
have precomputed a table of values of $f^{|C|^2 - |C| - 1}$, we can
see in time $O(n^3)$ whether the DNA Inequality holds in $\Gamma$.  So,
it suffices to show that we can compute the value of $f^{|C|^2 - |C| - 1}$
in time $O(n^{12}\log{n})$.  I claim that if we have a precomputed
table of values for $f^{k_1}$ and $f^{k_2}$, we can easily compute values
of $f^{k_1 + k_2}$, and can of course use this to precompute a table of
values for $f^{k_1 + k_2}$.

Say we wish to compute $f^{k_1 + k_2}(e_1, e_2)$. Consider the curve
with a first edge of $e_1$ and last edge of $e_2$, with  $k_1 + k_2 +2$
or fewer vertices, which has the minimal value of $f_\Gamma$. If we
consider an edge $e$ with at most $k_1 - 2$ vertices separating it from
$e_1$ and at most $k_2 - 2$ vertices separating it from $e_2$ (this
clearly exists as our curve has $k_1 + k_2 + 2$ or fewer vertices),
then we have that the value of $f_\Gamma$ of the entire curve is the
same as the sum of $f_\Gamma$ on the piece from $e_1 \to e$ plus the
value on the piece from $e \to e_2$, minus the length of $e$. Thus,
we have that
$f^{k_1 + k_2}(e_1, e_2) = \min_{e \in S}(f^{k_1}(e_1, e) + f^{k_2}(e,e_2) - |e|)$.
As $|S^3| = O(|C|^6) = O(n^{12})$, we have that to precompute a table
of values for $f^{k_1 + k_2}$ from a table of values for $f^{k_1}$ and
$f^{k_2}$ takes time $O(n^{12})$.  Using the double-and-add algorithm,
we can compute the table of values for $f^{|C|^2 - |C| - 1}$ in time
$O(n^{12}\log(|C|^2 - |C| - e)) = O(n^{12}\log{n})$, i.e.\ in polynomial
time.
\end{proof}

\begin{rem}
The above run-time analysis is quite pessimistic. By an earlier remark,
we actually have $|C| = O(n)$, so the above run-time analysis can be
improved to give $O(n^6\log{n})$.  We will leave the verification of
this claim to an interested reader.
\end{rem}

\section{\label{secseq} Sequences of Polygons}

\begin{define} \label{simplydented}
A polygon $\Gamma$ with a convex hull of $P$ is called \emph{simply dented} if,
for any two consecutive vertices of $\Gamma$, at least one is a vertex
of $P$, and for every two consecutive edges of $P$, at least one is an
edge of $\Gamma$.  Figure~\ref{figsimpdent} shows some examples of both
simply dented and non simply dented polygons.

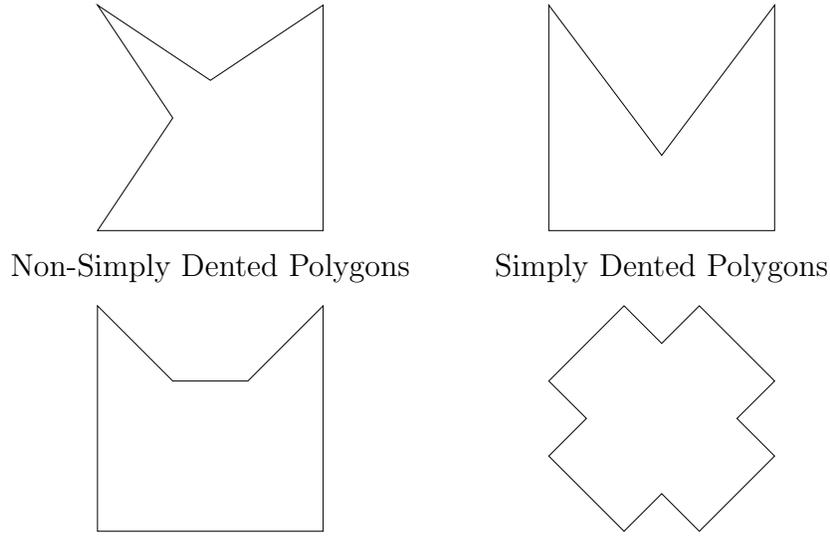
\begin{figure}[ht]
\begin{center}
\unitlength 1.0 mm
\begin{picture}(90,70)(0,0)
\path(0,0)(0,30)(10,20)(20,20)(30,30)(30,0)(0,0)
\path(0,40)(10,55)(0,70)(15,60)(30,70)(30,40)(0,40)
\path(60,40)(60,70)(75,50)(90,70)(90,40)(60,40)
\path(60,10)(65,15)(60,20)(70,30)(75,25)(80,30)(90,20)(85,15)(90,10)(80,0)(75,5)(70,0)(60,10)
\put(15,35){\makebox(0,0)[cc]{Non-Simply Dented Polygons}}
\put(75,35){\makebox(0,0)[cc]{Simply Dented Polygons}}
\end{picture}
\end{center}
\caption{Examples of Simply and Non-Simply Dented Polygons}
\label{figsimpdent}
\end{figure}
\end{define}

Fix some convex polygon $P$. Denote the set of points contained within
$P$ by $S$. Assume that all polygonal curves of interest are contained
within $P$, and have $M$ or fewer vertices. For every
$\mathbf{v} = (v_0, v_1, \ldots, v_{M-1}) \in S^M$, let
$\gamma[\mathbf{v}]$ be the closed curve $v_0v_1\cdots v_M=v_0$. Then,
any polygonal curve of interest is in the form $\gamma[\mathbf{v}]$
for some $\mathbf{v} \in S^M$. Note that the ``pseudo-vertices''
$v_i$ need not be real vertices of the curve $\gamma[\mathbf{v}]$,
as $v_{i-1},v_i,v_{i+1}$ might be collinear in that order for some
$i$. Moreover, consecutive $v_i$'s might be equal. Clearly, $S^M$ is a
compact space. It can be shown that:
\begin{itemize}
\item $\mathbf{v} \to \operatorname{perimeter}(\gamma[\mathbf{v}])$
is continuous on $S^M$.
\item $\mathbf{v} \to \operatorname{curvature}(\gamma[\mathbf{v}])$
is lower semicontinuous on $S^M$.
\end{itemize}
Now, consider what happens when we have a sequence of simply dented
polygons $P^1, P^2, P^3 \ldots$, with a common convex hull $P$.
Notice that there is some $M$ (twice the number of vertices of $P$
will do) such that each $P^k$ may be presented as $\gamma[\mathbf{V^k}]$
for $\mathbf{V} \in S^M$, since at least every other vertex of the $P^k$
is a vertex of $P$.

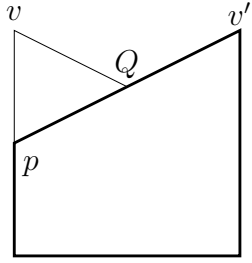
\begin{figure}[ht]
\begin{center}
\unitlength 1.5mm
\begin{picture}(20,27.5)(0,0)
\thinlines
\path(0,0)(0,20)(10,15)(20,20)(20,0)(0,0)
\put(0,21.5){\makebox(0,0)[cc]{$v$}}
\put(20,21.5){\makebox(0,0)[cc]{$v'$}}
\put(10,17){\makebox(0,0)[cc]{$Q$}}
\put(1.5,8){\makebox(0,0)[cc]{$p$}}
\Thicklines
\path(0,0)(0,10)(20,20)(20,0)(0,0)
\end{picture}
\end{center}
\caption{Critical Point Corresponding to $v'$}
\label{figc}
\end{figure}

Consider a vertex $v$ of $P^k$ which is also a vertex of $P$; it is
an endpoint of two edges of $P$. At most one of those edges of $P$ is
not an edge of $P^k$. If there is such an edge, we denote by $v'$ the
other endpoint of that edge.  If $Q$ is the vertex of $P^k$ between $v$
and $v'$, then we can form a non-vertex critical point by intersecting
$Qv$ with the boundary. We term this the critical point corresponding to
$v$. In Figure~\ref{figc} is pictured the critical point corresponding
to $v'$ (denoted $p$ in the diagram).

It is clear that when $k$ is sufficiently large, these are all of the
non-vertex critical points of $P^k$.

\begin{define} \label{gammakv}
We define the curve $\gamma^{k,v}$ to be the curve which is obtained
by starting with the curve $P^k$ and replacing the vertex $v$
with the critical point corresponding to $v'$, as pictured in
Figure~\ref{figc}. \end{define}

The next lemma essentially tells us that the DNA Inequality is true
for arbitrarily small dents of a region if it is true for these special
kinds of curves $\gamma^{k,v}$.

\begin{lm} \label{seq} Assume that we have a sequence of polygons
$P^1, P^2, P^3 \ldots$ such that:
\begin{itemize}
\item Each of the $P^k$ is simply dented;
\item None of the $P^k$ satisfy the DNA Inequality;
\item The $P^k$ have a common convex hull $P$;
\item For some $M$, there is a presentation $\mathbf{V^k} \in S^M$ of
each $P^k$ so that $\lim_{k \to \infty}{\mathbf{V^k}}$ exists in $S^M$,
and for which $\gamma[\lim_{k \to \infty}{\mathbf{V^k}}] = P$;
\item There exists $\epsilon > 0$, which does not depend upon $k$, such
that any two vertices of $P^k$ are at least $\epsilon$ apart for all $k$.
\end{itemize}
Then, it follows that there is an infinite subsequence of our sequence
in which there is a $CX_{P^k}$-polygon of the form $\gamma^{k,v}$, for
some $v$ which is both a vertex of $P$ (and of course consequently a
vertex of $P^k$), and is the endpoint of exactly one edge of $P$ which
is not also an edge of $P^k$ (for $k$ in our subsequence). \end{lm}

\begin{proof}
For $k$ sufficiently large, every critical point is either a vertex
or a critical point corresponding to the endpoints of some edge of the
convex hull not contained in $P^k$, since every two vertices of the $P^k$
are at least $\epsilon$ apart, and all of the $P^k$ are simply dented.
Thus, we will throw out the beginning of our sequence so that this is
true for all $k$.  For each $P^k$, we consider the set of edges of $P$
not contained in $P^k$. As there are finitely many possibilities for
this, there is an infinite subsequence such that the set is the same
for any element of the subsequence.

Thus it suffices to prove this lemma in the case where the set of edges
of $P$ not contained in $P^k$ does not depend on $k$.  We denote these
edges by $P_{r_\ell}P_{r_\ell+1}$ for $\ell = 1, 2, \ldots, \sigma$.  Now,
for sufficiently large $k$, $P^k$ is separable (because every two vertices
of the $P^k$ are at least $\epsilon$ apart), so we also assume that each
$P^k$ is separable. Then the vertex sequence of each $P^k$ is the same as
the vertex sequence for $P$ except that each edge $P_{r_\ell}P_{r_\ell+1}$
is replaced by $P_{r_\ell}Q^k_\ell P_{r_\ell+1}$ for some $Q^k_\ell$
in the interior of $P$.  Now, if we consider the elements
$(Q^k_1, Q^k_2, \ldots, Q^k_\sigma) \in S^\sigma$, and observe that
$S^\sigma$ is compact, it follows that we may select a subsequence in
which $(Q^k_1, Q^k_2, \ldots, Q^k_\sigma)$ converges in $S^\sigma$;
in other words, we may select a subsequence such that $Q^k_\ell$ has a
limit for each $\ell$.

Write $n$ for the number of vertices of $P$. To each $P^k$,
there is a $CX_\Gamma$-polygon $\gamma^k$. By the machinery of the
previous section, we may assume that $\gamma^k$ (as well as $P^k$)
have the number of vertices bounded by some function depending only
upon $n$, which we shall refer to as $M$.  By the remarks at the
beginning of the section, each of the $\gamma^k$ can be presented as
$\gamma[\mathbf{v^k}]$ for some $\mathbf{v^k} \in S^M$. (Recall that we
also notate $P^k$ as $\gamma[\mathbf{V^k}]$.)  As $S^M$ is compact, there
is a convergent subsequence of the $\mathbf{v^k}$, which converges to
$\mathbf{\overline{v}}$.  Let $AC$ represent average curvature, viewed as
a function from $S^M \to \mathbb{R}$. Since this is the product of a lower
semicontinuous function, curvature, and a continuous function, reciprocal
of perimeter, it's lower semicontinuous.  Now, since any two vertices of
$P^k$ are at least $\epsilon$ apart, it follows that $AC(\mathbf{V^k})$
converges to $AC(P)$.  Each $CX_\Gamma$-polygon $\gamma[\mathbf{v^k}]$
satisfies $AC(\mathbf{v^k}) < AC(\mathbf{V^k})$.  Since $\mathbf{v^k}
\to \mathbf{\overline{v}}$, and since $AC$ is lower semicontinuous,
$AC(\mathbf{\overline{v}}) \leq \limsup AC(\mathbf{v^k})$. Therefore:
$$AC(\mathbf{\overline{v}}) \leq \limsup AC(\mathbf{v^k}) \leq \limsup AC(\mathbf{V^k}) = AC(P).$$
Since $P$ is convex, the Lagarias-Richardson theorem \cite{L-R} tells us
that $AC(\mathbf{\overline{v}}) \geq AC(P)$, so $AC(\mathbf{\overline{v}}) = AC(P)$.
In \cite{N-P}, it is proven that for any convex curve $P$, the
only equality cases to the DNA Inequality are multiple circuits of $P$.
Therefore, $\gamma[\mathbf{\overline{v}}]$ is a multiple circuit of $P$.
We consider two cases:

\paragraph{Case 1:}
Our convergent subsequence contains infinitely many closed curves whose
vertex sequences contain non-critical points. By Lemma \ref{wmacritical},
we may assume that these have the form:
$$P^k_0, P^k_1, \ldots, P^k_{i-1}, P^k_i, X^k, P^k_j, P^k_{j+1}, \ldots, P^k_m = P^k_0$$
(which we will refer to as the first case) or
$$P^k_0, P^k_1, \ldots, P^k_{i-1}, P^k_i, X^k, Y^k, P^k_j, P^k_{j+1}, \ldots, P^k_m = P^k_0$$
(which we will refer to as the second case). From our convergent
subsequence, as there are finitely many choices for $i$, and $j$, we
may select a subsequence with $i$ and $j$ constant.

Now, I claim that the perimeter of $\gamma^k$ in this case is at most
the perimeter of $P$ plus twice the length of $P^k_iX^k$. This follows
from the triangle inequality: the length of $X^kP^k_j$ (respectively
$X^kY^k$) is less than or equal to the length of $P^k_iX^k$, plus the
length of the portion of the boundary between $P^k_i$ and $P^k_j$
(respectively $P^k_i$ and $Y^k$).  (Note that this argument relies
on $i \leq j$ in the first case or $i < j$ in the second case to
talk about the portion of the boundary between $P^k_i$ and $P^k_j$
or $P^k_i$ and $Y^k_j$.)  From this, it follows that the perimeter
of $\gamma[\mathbf{\overline{v}}]$ is at most the perimeter of
$P$ plus twice the diameter of $P$. Since twice the diameter of
$P$ is strictly less than the perimeter of $P$, the perimeter of
$\gamma[\mathbf{\overline{v}}]$ is strictly less than twice the perimeter
of $P$.  As $\gamma[\mathbf{\overline{v}}]$ is a multiple circuit of
$P$, it follows that $\gamma[\mathbf{\overline{v}}]$ is a single circuit
of $P$.  Since the perimeter of $\gamma[\mathbf{\overline{v}}]$ is the
same as the perimeter of $P$, we have:
$$\lim_{k \to \infty}(\mbox{perimeter of $\gamma^k$}) = \mbox{perimeter of $P$} = \lim_{k \to \infty}(\mbox{perimeter of $P^k$})$$
\begin{equation}\label{limcurv}
\Longrightarrow \lim_{k \to \infty}(\mbox{curvature of $\gamma^k$}) = \lim_{k \to \infty}(\mbox{curvature of $P^k$}) = 2\pi
\end{equation}

If $i = j$, elementary geometry shows that the (unsigned) curvature of
$\gamma^k$ is greater than $4\pi$. Thus, we may assume that $i < j$.
As $X^k$ is a turn-around, $X^k$ does not lie in the portion of $P^k$
between $P^k_i$ and $P^k_j$ (or $P^k_i$ and $Y^k$ in the second case).
Therefore, the only way for $\gamma[\mathbf{\overline{v}}]$ to be a
single circuit is for the length
$\dist(X^k, P^k_i) \to 0$ or $\dist(X^k, P^k_j) \to 0$.
(This should be replaced by $\dist(X^k, P^k_i) \to 0$ or
$\dist(X^k, Y^k) \to 0$ in the second case.)  If we are in the
second case, segment $X^kY^k$ intersects the boundary of $P^k$ in a
third point, say $Q^k_\ell$. Thus,
$\lim_{k \to \infty}(\dist(X^k,Y^k)) = \dist(P_{r_\ell},P_{r_\ell + 1})$.
From this we conclude that $X^k$ cannot approach $Y^k$. In other words,
we may assume without loss of generality that $\dist(X^k, P^k_i) \to 0$,
which implies that for $k$ sufficiently large, we have $X^k$ in either
$P^k_{i-1}P^k_i$ or $P^k_iP^k_{i+1}$.  If we are in the second case,
from elementary geometry it is clear that if $j = i + 1$ and $X^k$ is
between $Y^k$ and $P^k_j$ that we have the total curvature of $\gamma^k$
is greater than $4\pi$. Thus, by \eqref{limcurv}, we may assume that this
does not happen.  In either case, as $X^k$ does not lie in the portion
of $P^k$ between $P^k_i$ and $P^k_j$ (or $P^k_i$ and $Y^k$ in the second
case), we have that $X^k$ lies in the interval $P^k_{i-1}P^k_i$ for $k$
sufficiently large.  It follows that $\gamma^k$ has an angle with measure
$\pi$. (This occurs at $P^k_i$.)  As the limit of the total curvature
of $\gamma^k$ is $2\pi$, the limit of the sum of contributions to the
total curvature of every other angle is also $\pi$.  It follows that
$\gamma^k$ tends to some (degenerate) curve with two vertices, which is
not a multiple circuit of $P$.  Therefore, this case cannot happen.

\paragraph{Case 2:}
All but finitely many of the curves of our subsequence consist only of
critical points.  Throw out the beginning of our subsequence so that
all of the curves in the subsequence consist only of critical points.
Observe that the critical point corresponding to a vertex $v$ of $P$
(which of course is also a vertex of the $P^k$) tends to $v'$ as $k$
tends to $\infty$.  As $\gamma[\mathbf{\overline{v}}]$ is a multiple
circuit of $P$, it follows that for $k$ sufficiently large in our
subsequence, the vertices of the curve $\gamma^k$ are, in order,
(possibly for multiple circuits) exactly one of the (at most two;
one of them is $v$) critical points which becomes close to each vertex
$v$, and possibly visiting the $Q^k_\ell$ between $P^k_{r_\ell}$ and
$P^k_{r_\ell + 1}$.  For each vertex $v$ which is the endpoint of an
edge of $P$ which is not an edge of the $P^k$ in our subsequence,
write $n_{k,v}$ for the number of vertices of $\gamma^k$ which are
equal to the critical point corresponding to $v'$.  If for some $i$,
$\gamma^k_i$ is equal to the critical point corresponding to $v'$,
then I claim we may assume $v' \in \{\gamma^k_{i+1}, \gamma^k_{i-1}$\}.
For, if this is not the case, then we may replace $\gamma_i$ with $v$.
This increases the perimeter and leaves the curvature unchanged,
thus decreasing average curvature.  Therefore, we may assume that
when $\gamma^k_i$ equals any non-vertex critical point, then one
of $\gamma^k_{i \pm 1}$ is the vertex of $P$ to which the critical
point corresponds.  Suppose that $\gamma[\mathbf{\overline{v}}]$ is
a multiple circuit of $P$ which goes around $m$ times.  I claim that
$f_{P^k}(\gamma^k)=\sum n_{k,v}f_{P^k}(\gamma^{k, v})$, and this will
complete the proof since, as $f_{P^k}(\gamma^k) < 0$, it would follow
that one of the $f_{P^k}(\gamma^{k, v})$ is negative, for each $k$
in our subsequence that is sufficiently large.  To see the equality,
look at the two collections of curves:
\begin{itemize}
\item $n_{k,v}$ copies of $\gamma^{k,v}$ for each $k$ and $m$ copies of $P^k$.
\item $\sum{n_{k,v}}$ copies of $P^k$ and one copy of $\gamma^k$.
\end{itemize}
The sums over each collection of perimeter and of curvature are equal,
i.e., for each curve, compute the curvature and perimeter, then add
those values up. To see that the sums of the curvatures are equal,
look at the curvature contributions of the two collections at all the
possible vertices of the curves, and recall that when $\gamma^k_i$ equals
any non-vertex critical point, then one of $\gamma^k_{i \pm 1}$ is the
vertex of $P$ to which the critical point corresponds.  Similarly for
the perimeters, look at all possible edges of the curves.  By our earlier
comment, this completes the proof of this case, and of this lemma.
\end{proof}

\begin{rem}
If there is some way of verifying that the curves $\gamma^{k,v}$ are
not $CX_\Gamma$-polygons for large $k$, then the number of cases which
must be analyzed to directly apply this lemma in order to prove that the
DNA Inequality holds in $\Gamma$ is linear in the number of interior
vertices of $P$, which is significantly less than the number of cases
to directly apply Theorem \ref{poly}.
\end{rem}

\section{\label{secsmalldent} Classification of DNA-Polygons}
Here we prove our main result, Theorem~\ref{ismalldent} (Theorem~\ref{smalldent}).

\begin{thm} \label{smalldent}
If $P$ is a convex polygon with perimeter $p$ and we are denting an
edge with length $l$, and $\alpha$ is the larger of the two angles that
the edge makes with the two adjacent edges, then $P$ is a DDNA-polygon
(with respect to this edge) if and only if:
$$2p \leq \pi l \frac{1 + \cos{\alpha}}{\sin{\alpha}}.$$
\end{thm}

\begin{rem}
The set of convex polygons constructed in Theorem~\ref{smalldent}
is non-empty, as stated in the introduction. For example, it is easy
to see that it contains an isosceles right triangle (dented along the
hypotenuse). Of course, this proof is non-constructive, in the sense
that it does not tell by what angle you may dent a single edge. However,
by  Theorem \ref{poly}, we can compute what dents will work for any
specific curve. For example, for the isosceles right triangle, we can
find that for this case it holds as long as $\delta$ is less than or
equal to the root of
$$\frac{2\pi + 4\delta}{2 + \sqrt{2}\sec(\delta)} = \frac{4\pi + 6\delta}{4 - 2\tan(\pi/4 - \delta) + 2\sec(\pi/4 - \delta) + \sqrt{2}\sec(\delta)}$$
which is approximately $0.297142593$ radians.
\end{rem}

The following corollary gives a more verifiable way of testing if a
given convex polygon is a DDNA-polygon.

\begin{cor}
Let $P$ be a convex polygon with a fixed  edge and $\alpha$ be the larger
of the two angles that the fixed edge makes with the two adjacent edges.
Then, if $\alpha \leq \tan^{-1}(\pi/2) \approx 57.5^\circ$, $P$ is
a DDNA-polygon.  Additionally, if $P$ is a DDNA-polygon, then
$\alpha < \cos^{-1}\left(\frac{16 - \pi^2}{16 + \pi^2}\right) \approx 76.3^\circ$.
\end{cor}

\begin{proof}
To see the second statement, observe that
\begin{align*}
4l < 2p &\leq \pi l \frac{1 + \cos{\alpha}}{\sin{\alpha}} \\
\Longrightarrow 4\sin{\alpha} &< \pi (1 + \cos{\alpha}) \\
\Longrightarrow 16(1 - \cos^2{\alpha}) &< \pi^2 (1 + \cos{\alpha})^2 \\
\Longrightarrow \alpha &< \cos^{-1}\left(\frac{16 - \pi^2}{16 + \pi^2}\right) \approx 76.3^\circ.
\end{align*}

To see the first statement, consider such an $\alpha$ and the isosceles
triangle with equal angles $\alpha$, and base $l$. Then, the perimeter $p$
of $P$ is less than or equal to the perimeter of our isosceles triangle,
which equals $l\left(1 + \frac{1}{\cos{\alpha}}\right)$.  So it suffices
to have
$$2l\left(1 + \frac{1}{\cos{\alpha}}\right) \leq \pi l \frac{1 + \cos{\alpha}}{\sin{\alpha}},$$
which is equivalent to $\alpha \leq \tan^{-1}(\pi/2) \approx 57.5^\circ$.
\end{proof}

\noindent
Now, we prove Theorem~\ref{smalldent}.

\begin{proof}
Recall Definition~\ref{Pdelta} and the notation there. Then $P$ is \emph{not}
a DDNA-polygon if and only if there exists a sequence
$\delta_1, \delta_2, \delta_3, \ldots$ with a limit of 0 such
that there is a $CX_{P_{\delta_k}}$-polygon for each $\delta_k$.
Denote the edge that we are denting by $AB$.  By Lemma~\ref{seq}, we
have that this happens if and only if, for $\delta$ arbitrarily small,
one of the curves $\gamma^{k, A}, \gamma^{k, B}$ (which are pictured in
Figure~\ref{figgkab}) is a $CX_{P_{\delta_k}}$-polygon.

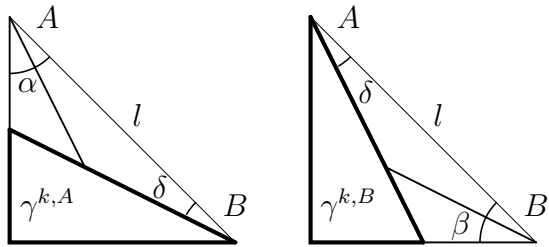
\begin{figure}[ht]
\begin{center}
\unitlength 1.0mm
\begin{picture}(70,30)(0,0)
\multiput(5,30)(40,0){2}{\makebox(0,0)[cc]{$A$}}
\multiput(30,5)(40,0){2}{\makebox(0,0)[cc]{$B$}}
\multiput(17,17)(40,0){2}{\makebox(0,0)[cc]{$l$}}
\multiput(20,7.4)(27.4,12.6){2}{\makebox(0,0)[cc]{$\delta$}}
\put(2.4,21){\makebox(0,0)[cc]{$\alpha$}}
\put(60,2.25){\makebox(0,0)[cc]{$\beta$}}
\put(5,5){\makebox(0,0)[cc]{$\gamma^{k, A}$}}
\put(45,5){\makebox(0,0)[cc]{$\gamma^{k, B}$}}
\allinethickness{0.15mm}
\path(0,30)(30,0)
\path(40,30)(70,0)
\allinethickness{0.25mm}
\put(40,30){\arc{15}{0.7854}{1.0472}}
\put(30,0){\arc{15}{3.6652}{3.927}}
\put(0,30){\arc{15}{0.7854}{1.5708}}
\put(70,0){\arc{15}{3.1416}{3.927}}
\allinethickness{0.3mm}
\path(0,15)(0,30)
\path(0,30)(10,10)
\path(50,10)(70,0)
\path(70,0)(55,0)
\allinethickness{0.6mm}
\path(0,15)(0,0)(30,0)
\path(40,30)(40,0)(55,0)
\path(0,15)(30,0)
\path(55,0)(40,30)
\end{picture}
\end{center}
\caption{The Curves $\gamma^{k, A}$ and $\gamma^{k, B}$}
\label{figgkab}
\end{figure}

If we write $l$ for the length $AB$, and $\alpha, \beta$ for the angles
at $A, B$ in $P$, then $f_{P^k}(\gamma^{k, A}) < 0$ if and only if:
$$\frac{2\pi}{p + l\frac{\sin(\alpha)}{\sin(\alpha + \delta)} - l - l\frac{\sin(\delta)}{\sin(\alpha + \delta)}} < \frac{2\pi + 4\delta}{p + l(\sec(\delta) - 1)},$$
where $p$ is the perimeter of $P$.  Similarly, 
$f_{P^k}(\gamma^{k, B}) < 0$ if and only if the above is true with
$\alpha$ replaced by $\beta$. Therefore, the DNA Inequality holds for
arbitrarily small dents if and only if, for $\delta$ arbitrarily small,
we have:
$$\frac{2\pi}{p + l\frac{\sin(\alpha)}{\sin(\alpha + \delta)} - l - l\frac{\sin(\delta)}{\sin(\alpha + \delta)}} \geq \frac{2\pi + 4\delta}{p + l(\sec(\delta) - 1)}$$
where $\alpha$ assumes either angle.  Now, as $LHS(0) = RHS(0)$, the
above holds for $\delta$ arbitrarily small implies that
$LHS'(0) \geq RHS'(0)$; if $LHS'(0) > RHS'(0)$, then the above holds for
arbitrarily small $\delta$.  We compute first and second derivatives at 0:
\begin{align*}
LHS'(0) &= \frac{2\pi l (1 + \cos{\alpha})}{p^2\sin{\alpha}} \\
LHS''(0) &= \frac{2\pi l (1 + \cos{\alpha})^2 (2l - p)}{p^3\sin^2{\alpha}} \\
RHS'(0) &= \frac{4}{p} \\
RHS''(0) &= \frac{-2\pi l}{p^2}
\end{align*}
Thus, $P$ is a DDNA-polygon if:
\begin{align*}
\frac{2\pi l (1 + \cos{\alpha})}{p^2\sin{\alpha}} &> \frac{4}{p} \\
\Longleftrightarrow 2p &< \pi l \frac{1 + \cos{\alpha}}{\sin{\alpha}}
\end{align*}
and only if $2p \leq \pi l \frac{1 + \cos{\alpha}}{\sin{\alpha}}$.
If we observe that $\frac{1 + \cos{\alpha}}{\sin{\alpha}}$ is a decreasing
function (its derivative is $-\frac{1+\cos{\alpha}}{\sin^2{\alpha}}$),
it follows that we may assume that $\alpha$ is the bigger of the two
angles for the above two statements. From here on, we assume this.

I claim that in the equality case, the DNA Inequality holds. As
$LHS(0) = RHS(0)$ and $LHS'(0) = RHS'(0)$, it suffices to
examine the second derivative.  Assume that
$2p = \pi l \frac{1 + \cos{\alpha}}{\sin{\alpha}}$. The fact that
$l(1 + \sec{\alpha}) \geq p$ (which comes from the fact that $P$ is
contained in an isosceles triangle with base $l$ and angles $\alpha$
at the base) implies that
$2l(1 + \sec{\alpha}) \geq \pi l \frac{1 + \cos{\alpha}}{\sin{\alpha}}$,
from which it follows that $\tan{\alpha} \geq \pi/2$, with
$l(1 + \sec{\alpha}) = p$ if and only if $P$ is an isosceles triangle with
base $l$, and angles $\alpha$ at the base. A simple calculation, using
$LHS''(0)$ and $RHS''(0)$ above shows that $LHS''(0) > RHS''(0)$ if and
only if $\tan{\alpha} > \pi/2$. Thus, it suffices to examine the case of
an isosceles triangle with base $4$ and height $\pi$.  For this triangle,
we can explicitly compute:
\begin{align*}
LHS &= \frac{2\sin{\delta} + \pi\cos{\delta}}{2 + \sqrt{4 + \pi^2}\cos{\delta}} \\
RHS &= \frac{(\pi + 2\delta)\cos{\delta}}{2 + \sqrt{4 + \pi^2}\cos{\delta}}
\end{align*}
and find that $LHS > RHS$, for all $\delta$ for $\pi/2 > \delta > 0$.
\end{proof}

\section{\label{secclos} Concluding Remarks}

In this paper, we have assumed for simplicity that the dent $XAB$ in
definition~\ref{Pdelta} is isosceles. However, as noted in Section~\ref{secintro}, the
methods of this paper can still be applied if the triangle is not
isosceles, or even if there are multiple dents all depending on one
parameter $\delta$, so long as adjacent sides are not dented and there
exists $\epsilon > 0$, which does not depend on $\delta$, such that any
two vertices of $P_\delta$ are at least $\epsilon$ apart. Essentially,
so long as a sequence of counterexamples $P_{\delta_k}$ would meet the criteria of Lemma~\ref{seq},
we can analyze the DNA Inequality in $P_\delta$ using the methods of this paper.

It is reasonable to conjecture that Theorem~\ref{smalldent} holds
in cases where the curve is not polygonal, even though the techniques
of this paper can probably not be used to prove it.
More precisely, let $P$ be a piecewise-smooth convex curve, with at
least one of the pieces straight (call this piece $AB$).  Write $p$
for the arc length of $P$, and $\alpha$ for the larger of the two angles
that $AB$ makes with the other one-sided tangent vector to $A$ and $B$.
(For example, if $P$ is smooth, $\alpha = \pi$. If $P$ is a semicircle
and $AB$ is the diameter, then $\alpha = \pi/2$.)  Then, based on
Theorem~\ref{smalldent}, we conjecture that the DNA inequality holds
for arbitrarily small dents of $P$ along side $AB$ if and only if
$$2p \leq \pi l \frac{1 + \cos{\alpha}}{\sin{\alpha}},$$
where $l$ is the length of $AB$.

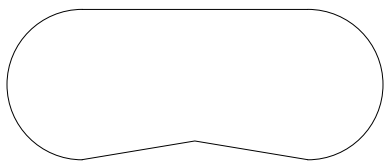
\begin{figure}[ht]
\begin{center}
\unitlength 1 mm
\begin{picture}(50,20)(0,0)
\path(10,20)(40,20)
\path(10,0)(25,2.5)(40,0)
\put(10,10){\arc{20}{1.5708}{4.7124}}
\put(40,10){\arc{20}{4.7124}{7.854}}
\end{picture}
\end{center}
\caption{A Dent in a Smooth Curve}
\label{stadium}
\end{figure}

In particular, if true, this conjecture would imply that the DNA
Inequality fails for arbitrarily small dents of smooth convex curves.
(An example is pictured in Figure~\ref{stadium}.)  While the above
conjecture is probably quite difficult to prove, this corollary for
arbitrarily small dents of smooth convex curves is probably not too
difficult to prove: the correct curve to take for the counterexample
DNA should be analogous to the $\gamma^{k, A}$ used in the proof of
Theorem~\ref{smalldent}.

\end{document}